\documentclass{article}

\usepackage[margin=2.54cm]{geometry}

\usepackage[english]{babel}
\usepackage{amsmath,amssymb,amsthm}
\usepackage{bbm}
\usepackage{parskip}
\usepackage{tikz}
\usetikzlibrary{automata, arrows.meta, positioning}
\usepackage{xcolor}
\usepackage{cite}
\usepackage{subcaption}
\usepackage{bbm}
\usepackage{rotating}
\usepackage{mathtools}
\usepackage{graphicx}
\graphicspath{ {./images/} }

\usepackage{hyperref}

\newtheorem{theorem}{Theorem}[section]
\newtheorem{corollary}{Corollary}[theorem]
\newtheorem{lemma}[theorem]{Lemma}

\newtheorem{prop}{Proposition}[section]
\newtheorem*{theorem*}{Theorem}

\theoremstyle{definition}
\newtheorem{definition}{Definition}[section]
\newtheorem{conjecture}{Conjecture}[section]

\newtheorem*{definition*}{Definition}

\newcommand{\N}{\mathbb{N}}
\newcommand{\Z}{\mathbb{Z}}
\newcommand{\Q}{\mathbb{Q}}
\newcommand{\R}{\mathbb{R}}

\newcommand{\D}{\mathcal{D}}

\DeclareMathOperator{\var}{var}
\DeclareMathOperator{\pos}{pos}
\DeclareMathOperator{\zero}{zero}
\DeclareMathOperator{\dist}{dist}

\newcommand*{\defeq}{\mathrel{\vcenter{\baselineskip0.5ex \lineskiplimit0pt
                        \hbox{\scriptsize.}\hbox{\scriptsize.}}}%
        =}

\title{Slicing the torus, thermodynamic formalism of self similar measures with overlaps}
\author{Peej Ingarfield}
\date{01/04/2026}

\begin{document}
\maketitle

\begin{abstract}
Iterated function systems of the form \(\{\frac{x+i}{\beta} : i \in \{0,1,f_{1},\dots,f_{k}\}\in \Q[\beta]\}\) give a large and natural class of self similar measures in \(\R\) that exhibit overlaps. 
In the case that \(\beta \in (1,2]\) is a  Pisot–Vijayaraghavan number, \(k>1\) and \(f_{i}\in \Q[\beta]\), this family of such measures exhibit exact overlap and are a generalisation of the well studied Bernoulli convolutions. 
We show that this family of generalised Bernoulli convolution exhibit dimension drop.
The key result of this work relates a pressure function on the torus to an upper bound for the amount dimension can drop for measure in this family. 
Of note is how varying the parameters \(f_{1},\dots,f_{k}\) changes the rational slice the pressure function is restricted to but does not change the potential function used. 
\end{abstract}

\section{Introduction}
\begin{figure}[h!]
   \centering
    \begin{subfigure}[t]{0.5\textwidth}
        \centering
        \includegraphics[height=4.5cm, trim={4cm 10cm 4cm 10cm},clip]{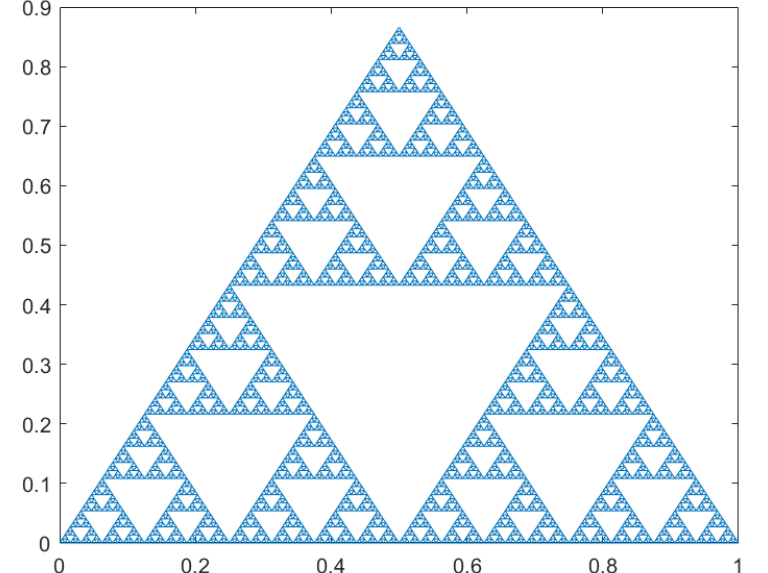}
        \caption{The Sierpinski Triangle\label{fig:sierpinski}}
    \end{subfigure}%
    ~ 
    \begin{subfigure}[t]{0.5\textwidth}
        \centering
        \includegraphics[height=4.5cm,trim={4cm 10cm 4cm 10cm},clip]{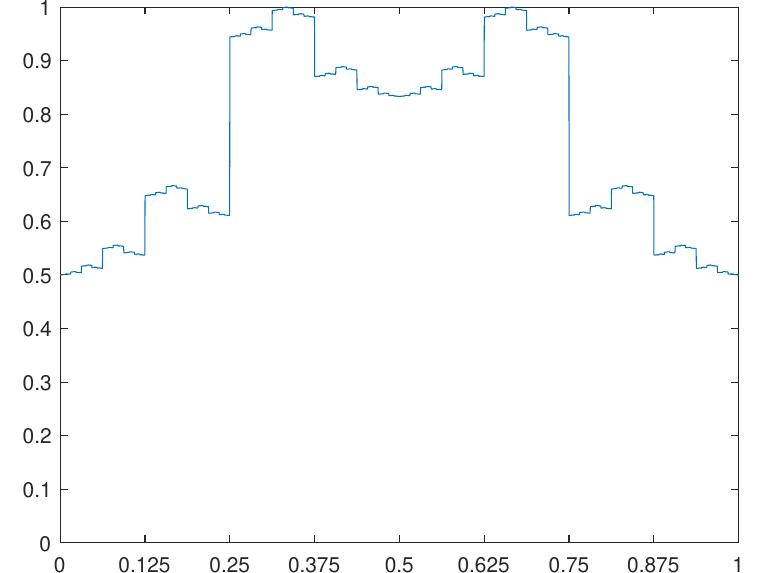}
        \caption{Radon–Nikodym Derivative of the vertical projection of the Sierpi\'nski Triangle into \(\R\).}
    \end{subfigure}
\end{figure}

Suppose that the functions $g_1,\dots,g_k\colon \R^d\to\R^d$ are contracting similarities 
i.e., \(\lVert g_j(x) - g_j(y)\rVert = r_j \lVert x-y\rVert\) for all \(x,y\in\R^d\), where \(r_j\in (0,1)\). 
Then \(\{g_1,\dots,g_k\}\) is called an \emph{iterated function system} (IFS) of similarities. 
For such an IFS, there is a unique non-empty compact \emph{attractor}. 
This is the set
 \(X\subset \R^d\), satisfying 
   \[ X = \bigcup_{j=1}^k g_k(X).\]
   (See~\cite{hutchinson1981fractals}.) 
   Many classical objects arise as attractors of such systems. 
   For example, the Sierpi{\'n}ski triangle (Figure~\ref{fig:sierpinski}) is the
   attractor of the IFS given by \(\{g_0,g_1,g_2\}\), where \(g_i(x)=x/2 + t_i\), 
   and \(t_0=(0,0)\), \(t_1=(1/2,0)\) and \(t_2 = (1/4,\sqrt{3}/4).\)

Separation conditions have played a key role in the understanding of such fractal attractors. 
Two common separation conditions are the \emph{strong separation condition} (SSC) and the
weaker \emph{open set condition} (OSC). 
The OSC requires that there exists an 
non-empty open set \(U\subset\R^d\) such that the images \(g_{i}(U)\) are pairwise disjoint. 
For the Sierpi{\'n}ski triangle, such \(U\) is given by an equilateral triangle with vertices
\((0,0)\), \((1,0)\) and \((1/2,\sqrt{3}/2).\)

If the OSC holds, then Hutchinson's formula \cite{hutchinson1981fractals} gives a simple expression for the attractor's Hausdorff dimension. 
On the other hand, if an IFS does not satisfy the open set condition, we may consider it to have \emph{overlaps}. 
Overlaps in fractals are typically considered to be one of two types. 
For an IFS of similarities \(\{g_{1},\dots,g_{k}\}\), the IFS is said to have \emph{exact overlaps} if there exists a \(n\in \N\) such that \(g_{i_{1}}\circ\dots\circ g_{i_{n}} = g_{j_{1}}\circ\dots\circ g_{j_{n}}\) where \(i_{1}\dots i_{n} \neq j_{1}\dots j_{n}.\)
If an IFS does not meet a separation criteria nor have exact overlaps it is said to have \emph{non exact overlaps}. 

Some IFSs with overlaps may lead to attractors with uninteresting or well understood geometry. 
The IFS \(\{x\to\frac{x+i}{2}: i \in\{0,1,1/2\}\}\) has the attractor of \([0,1]\), the interval has dimension 1 and is well understood.  
Consequently rather than consider attractors we study a probability measure supported on the attractor. 
For the IFS of similarities \(\{g_{1},\dots,g_{k}\}\) with attractor \(X\) we study the unique probability measure satisfying ,
\[\mu(X) = \frac{1}{k}\sum_{i=1}^{k}\mu(g^{-1}_{i}(X)).\]

It is well known that there are many different kinds of dimension for fractal attractors. 
For measures the same is true. 
Rather than introduce many of these different notions for dimension we shall appeal to general theory to motivate working with a given definition of Hausdorff dimension.

Work by Feng and Hu \cite{Feng} gives that the measure, \(\mu\) for an IFS of similarities is \emph{exact dimensional} , which means \(\lim_{r\to0} \frac{\mu(B(x,r))}{r}=\alpha\) for \(\mu\) almost every \(x\). 
While the work of Fan, Lau, Rao \cite{measuredim} give a rich theory of relations for different dimensions of a measure. 
Notably for our context, should a measure be exact dimensional taking value \(\alpha\) then its Hausdorff dimension exists and takes value \(\alpha\).
We define the \emph{Hausdorff dimension of a measure} as \[\dim(\mu) \defeq \inf\{\dim(E) \colon \mu(E)>0\}.\]
This would generally be the definition for \emph{lower Hausdorff dimension } however by the above the Hausdorff dimension exists and so lower Hausdorff dimension and Hausdorff dimension take the same value. 

Exact overlaps are of note due to the key role they play in the dimension of a measure. 
This can be seen in a folklore conjecture known as the exact overlap conjecture. 
To concisely state this conjecture we must first introduce entropy and Lyapunov exponent. 
Let \(G = \{g_{1},\dots,g_{k}\} \) be the IFS of similarities chosen with positive probabilities \(\mathbf{p}\).
Then the \emph{Lyapunov exponent} ,\(\lambda_{G}(\mathbf{p})\), and \emph{entropy}, \(H_{G}(\mathbf{p})\), of the IFS \(G\) with probabilities \(\mathbf{p}\) are defined as \[\lambda(\mathbf{p}) \defeq-\sum_{j=0}^{k}p_{j}\log r_{j},\quad H(\mathbf{p}) \defeq -\sum_{j=0}^{k}p_{j}\log p_{j}.\]
The \(G\) in the subscript will often be suppressed when it is clear from context what IFS is being considered. 

\begin{conjecture}\label{conj:exact overlaps}
    For an IFS of similarities of \(\R^{1}\)
    then \(\dim(\mu) = \min\{1,\frac{H(\mathbf{p})}{\lambda({\mathbf{p})}}\}\) unless the IFS has exact overlaps.
\end{conjecture}

In the case that \(\dim(\mu) < \min\{1,\frac{H(\mathbf{p})}{\lambda({\mathbf{p})}}\}\)  then we say that the measure \(\mu\) or the IFS has \emph{dimension drop}. For work on the exact overlap conjecture from a perspective focusing on when dimension does not drop see Rapaport,\cite{rapaport2020proofexactoverlapsconjecture}.

Calculation of the Hausdorff dimension is in general a difficult problem, whether studying measures with or without overlaps. 
The \emph{\(L^{q}\) dimensions} are a one parameter family of values that are useful in categorising how smooth a measure is at varying scales and, in bounding Hausdorff dimension. 
To understand how \(L^{q}\) dimension does this we introduce the objects to define this family. 
Let \(\Lambda\) be a partition which generates the Borel sigma algebra of the \([0,\frac{1}{\beta-1})\) for \(\beta \in (1,2]\) for the map \(\sigma\).
Then let \(\Lambda_{n}\) be the \(n^{th}\) iterate of generating the sigma algebra,\(\Lambda_{n} \defeq \bigvee\sigma^{-n}(\Lambda). \)

From this we can define the \(L^{q}\) dimension of the \(\sigma\) invariant measure \(\mu\) as,
\[D_{\mu}(q)\defeq\liminf_{n\to\infty}\frac{-\log \sum_{I\in \Lambda_{n}}\mu(I)^{q}}{n(q-1)}.\]

The sum \(\sum_{I\in \Lambda_{n}}\mu(I)^{q}\) can be thought of as quantifying smoothness as measure.
This is as the sum takes a maximal value for a uniform measure and decreasing values as the measure concentrates mass toward points. 

In the case that the measure \(\mu\) is self similar it is known from the work of Shmerkin, \cite{pabloselfsimlq} Theorem 5.1, that \(D_{\mu}(q)\to\dim(\mu)\) from below as \(q\to 1\) from above. 

Work by Lau, Ngai and Rao \cite{laulq}, gives a computable algorithm for calculation of the \(L^{2}\) dimension of a wide class of self similar measures with overlaps. 
The measure we will consider are a sub class of those considered in their work.
Their algorithm sensitively depends on the parameters of the system. 
Changing any of the values of the system involved means recomputing a large bespoke matrix. 
We shall deviate from the approach of \(L^{2}\) dimension and construct a bound which is always lower than the \(L^{2}\) dimension but has a less sensitive dependence on the system's parameters.

Before introducing the family we shall consider we provide a brief overview of motivating results for Bernoulli convolutions, a family of measures with exact overlaps.
The family we will study can be considered as a generalisation of Bernoulli convolutions to a higher number of base masses.

Let \(\mathcal{D}\) be a set of at least two elements, \(\beta > 1\) a fixed constant and \(\nu_{\beta,\mathcal{D}}\) be the weak star limit of the following sum of Dirac masses
\begin{equation}\label{equation : weakstar}
    \nu_{\beta,\mathcal{D}} = \lim_{n\to\infty} \nu_{\beta,\mathcal{D}}^{(n)} = \lim_{n\to\infty} \frac{1}{\lvert\mathcal{D}\rvert^{n}}\sum_{a_{1}\dots a_{n}=a \in\mathcal{D}^{n}}\delta_{\sum_{i=1}^{n}a_{i}\beta^{-i}}.
\end{equation}

The measures \(\nu_{\beta,\mathcal{D}}\) are equicontractive self similar measures of the line.
These measures, \(\nu_{\beta,\mathcal{D}}\) can be seen as the measures associated to the IFS formed by \(\{x\to\frac{x+i}{\beta} :  i\in\mathcal{D}\}\). 
In the case that \(\mathcal{D} = \{0,1\}, \beta \in (1,2)\) these equicontractive self similar measures are a well studied family of fractal measures known as \emph{Bernoulli convolutions}.
Bernoulli convolutions are self similar measure with overlaps, the study of such measures is an active topic of research.
A central question is whether a given Bernoulli convolution is absolutely continuous or not with respect to the Lebesgue measure. 
If a given measure can be shown to be singular then it is often asked whether the measure has dimension \(< 1.\)

It is known that the algebraic properties of \(\beta\) are key to the dimension of Bernoulli convolutions. 
An algebraic number is called a \emph{ Pisot–Vijayaraghavan number}, or a \emph{PV number}, if it is a real algebraic integer greater than 1, whose Galois conjugates are all less than 1 in absolute value. 
The minimum polynomial of such a \(\beta\) is a monic integer polynomial and is often denoted \(p(\beta).\) 
The order of \(\beta\), denoted \(m(\beta)\), is the degree of \(p(\beta).\)

Erd\"os \cite{Erdos} showed that the Bernoulli convolution \(\nu_{\beta,\{0,1\}}\) is singular when \(\beta\) is a Pisot–Vijayaraghavan number. 
This result was furthered by Garsia \cite{garsia} who showed that the Hausdorff dimension of a Bernoulli convolution is less than one for such \(\beta\). 
So far these are the only known examples of Bernoulli convolutions of dimension less than one.
Garsia in \cite{Garsia2} also constructed explicit examples of a family of \(\beta\) for which \(\nu_{\beta,\{0,1\}}\) are absolutely continuous. 
Beyond explicit examples we note the work of Solomyak \cite{Solomyak} which showed that Bernoulli convolutions are absolutely continuous for almost all \(\beta\in(0,1)\). 
This was followed by the work of Shmerkin \cite{shmerkin} which proves that Bernoulli convolutions are absolutely continuous except on a set of dimension zero. 
Breuillard and Varju \cite{Breuillard} gave an explicit uniform lower bound on the dimension of Bernoulli convolutions for all algebraic integers, \(\beta \in(1,2)\). 

We now make a change of notation, for the sake of legibility in avoiding overly cluttered subscripts. 
Let \(\Q[\beta]\) be the field formed by adjoining \(\beta\) to \(\Q\).
For \(f_{1},\dots,f_{k}\in \Q[\beta]\), let \(\underline{\theta} = (f_{1},\dots,f_{k})\).
Then \[\mu_{\beta,\underline{\theta}}\defeq \nu_{\beta,\{0,1,f_{1},\dots,f_{k}\}}.\]
We shall study the family of unique probability measures that satisfy the following,

\begin{equation}\label{equation:mu self sim}
\mu_{\beta,\underline{\theta}}(A) = \frac{1}{k+2}\left( \mu_{\beta,\underline{\theta}}(\beta A) + \mu_{\beta,\underline{\theta}}(\beta A - 1) +\sum_{i=1}^{k} \mu_{\beta,\underline{\theta}}(\beta A - f_{i}) \right).
\end{equation}

An atypical example of the measures we shall considering can be seen as projections of the natural measure on the Sierpi\'nski triangle. This measure was seen as an example in (\!\!\cite{Hochman_inverse} pg 5). 
This example is useful for providing a good intuition and visualisation of the family of measures but its behaviour is unlike the rest of the family of measures.
The behavioural difference is due to its integer contraction ratio instead of the algebraic integer contraction resulting in the digits \(0,1\) not having overlap in base \(2\) while having overlap for all algebraic integers considered.

As we are considering contractions with ratio \(\beta\) we introduce the following notion akin to the doubling, the \emph{\(\beta\)ing map} , \(T_{\beta}:\mathbb{T}^{k}\to\mathbb{T}^{k}\) by \(T_{\beta}(x) = \beta x.\)

Motivated by the above results, we show in section \ref{section : dim drop} that this generalisation of the Bernoulli convolution to the family \(\mu_{\beta,\underline{\theta}}\) exhibits dimension drop for a PV contraction ratio and rational digits. 
To the best of the authors knowledge this result is new to the literature, but can be readily derived from a combination of more general results.

\begin{theorem}
    For \(\beta\in(1,2]\) a Pisot Vijayaraghavan number, \(f_{1},\dots,f_{k}\in \Q[\beta]\) ,\(k \geq 1\) and \(\underline{\theta} = (f_{1},\dots,f_{k})\) then \[\dim(\mu_{\beta,\underline{\theta}})<1.\]
\end{theorem}

This leads to the natural question of how much the dimension can drop?
We give an upper bound for the amount the dimension can drop in terms of the pressure function of a specified dynamical system and potential function. 
This gives a new way of understanding dimension drop in terms of varying subspaces of the torus. 
Of note is that the potential function exists on the \(k+1\) torus and is independent of the choice of \(\beta,f_{i},\dots,f_{k}\). 
This can be seen in our main theorem, Theorem \ref{theorem: dim drop bound}, which we state a version of below.

\begin{theorem}
    There exists \(\phi :[-\frac{1}{\beta-1},\frac{1}{\beta-1})^{k} \to \R\), such that for every \(f_{i},\dots f_{k} \in \Q[\beta]\) and \(\underline{\theta} = (f_{1},\dots,f_{k})\), \[1 > \dim(\mu_{\beta,\underline{\theta}} )\geq \frac{2\log(k+2)}{\log\beta} -\frac{P(l_{\underline{\theta}},T\vert_{l_{\underline{\theta}}},\phi)}{\log \beta}>0, \]
    where \(T\vert_{l_{\underline{\theta}}}\) is the \(\beta\)ing map restricted to the plane at rational angle \(\underline{\theta}\) on the k+1 torus, \(\R^{k+1}/\Z^{k+1}\) , and \(P(l_{\underline{\theta}},T\vert_{l_{\underline{\theta}}},\phi)\) is the topological pressure of \(\phi\) under the map \(T\vert_{l_{\underline{\theta}}}\).
\end{theorem}

We now provide a brief outline of the proof of this theorem. 
A key result in Hochman's work on measures with exact overlaps \cite{Hochman} is that the dimension of the measures can be expressed as the minimum of 1 and the ratio of the random walk entropy and Lyapunov exponent. 
For the definition of random walk entropy see definition \ref{definition : random walk entropy}.
Rather than work with random walk entropy we pass to considering the growth rate of exact overlaps.
This is done by modifying the arguments used in the work of Akiyama, Feng, Persson and Kempton (\!\!\cite{Akiyama} prop 3.5).
From this, we are able to carefully construct a potential function, \(\phi\), on \([-\frac{1}{\beta-1},\frac{1}{\beta-1})^{k+1}\).
This \(\phi\) counts the growth rate of the number of exact overlaps in \(\mu_{\beta,\underline{\theta}}\) when restricted to an appropriate subspace of the \(k+1\) torus.
The construction of the final dynamical system and ensuring the potential function has the required properties is the bulk of the work of this paper.  

There are many methods to calculate dimension that rely on potential functions such as \cite{barany2014dimension}.
The potential function we give in definition \ref{definition : potential}, has some key differences to other methods in the literature.
The matrices that the potential function is defined by do not change as the parameters \(\beta,f_{1},\dots,f_{k}\) are varied. 
Rather than change the matrices, we vary the subspaces of the torus which we restrict our pressure function to.
This means that by gaining a deeper understanding of how products of a fixed number of matrices behave, we may understand bounds for dimension drop. 
We hope this can provide a new avenue to understanding dimension drop of more complex systems and give a new method for estimating dimension drop for existing ones.

Finally in section \ref{gibbsprop}, we explore the Gibbs properties of the measure \(\mu_{\beta,\underline{\theta}}.\)
Due to the technical nature of the definitions involved rather than state this result here we direct the reader to section \ref{gibbsprop} for the details.

\section{Preliminaries \& Notation}
\subsection{Symbolic Dynamics}
For a given finite alphabet, \(A \defeq \{a_{1},\dots,a_{l}\}\), we denote the space of all infinite sequences over \(A\) as \(A^{\N}\). 
Further, denote the space of all finite words over \(A\) as \(A^{*}\) and the space of all words of length exactly \(n\) for each \(n\in\mathbb{N}\) as \(A^{n}\). 
We denote the \(i^{th}\) letter of a word, \(\omega\), in any of these spaces as \(\omega_{i}\).
For \(\underline{\omega} \in A^{\N}\), \( \underline{\omega} = \omega_{1} \omega_{2}  \dots \) then let \(\sigma\) be the left shift defined by \( \sigma( \underline{\omega} ) \defeq \omega_{2} \omega_{3} \dots \). 
Similarly, we define \(\sigma(\omega) \defeq \omega_{2}\dots\omega_{n}\) in the case of finite words. 
For a word \(\omega_{1}\dots\omega_{n} \in A^{*}\) define the cylinder set as \([\omega_{1}\dots\omega_{n}] \defeq \{ \gamma \in A^{*} 
\colon \omega_{1}\dots\omega_{n} = \gamma_{1}\dots\gamma_{n} \}\). 
We define cylinder sets for \(A^{\N}\) analogously. 

These symbolic spaces are given geometric meaning through the use of the the projection map, \(\pi_{\beta} \colon A^{\N} \to \R\) by \[\pi_{\beta}(\underline{\omega}) \defeq \sum_{i=1}^{\infty} \beta^{-i}\underline{\omega}_{i}\]
For \(n\in \N\) we define \(A^{n}\) analogously with a finite sum. 
At times we may drop the subscript \(\beta\), this is done to aid readability. 
In the case that \(\beta\) is suppressed we are considering a fixed but arbitrary PV number such that \(1<\beta\leq2\). 

Much of the dynamics we study are on cubes of various dimensions, to which we associate symbolic spaces.
Let \(B_{k} \defeq \left\{ (i_{1},\dots,i_{k}) \colon i_{j} \in \{-1,0,1\} \forall j \right\}, \) this can be seen as a k dimensional version of the signed binary alphabet.
We define the \emph{k dimensional \(\beta\) expansion} of a point \((x_{1},\dots,x_{k}) \in [-\frac{1}{\beta-1},\frac{1}{\beta-1})^{k}\) as the sequence \(\underline{\omega} \in B_{k}^{\N}\) such that \(\pi(\underline{\omega}) = \sum_{i} \omega_{i}\beta^{-i} = (x_{1},\dots,x_{k})\). 
The space \([-\frac{1}{\beta-1},\frac{1}{\beta-1})^{k}\) is equipped with the \(\beta\)ing map \(T_{\beta} : [-\frac{1}{\beta-1},\frac{1}{\beta-1})^{k} \to [-\frac{1}{\beta-1},\frac{1}{\beta-1})^{k}\), defined as
\[T_{\beta}(x_{1},\dots,x_{k}) = (\beta x_{1} \mod \frac{1}{\beta-1} ,\dots, \beta x_{k}\mod \frac{1}{\beta-1}).\]
Let \( \mathbb{T}^{k} \defeq \mathbb{R}^{k} / \mathbb{Z}^{k}\) be the usual k-torus which we identify with \([-\frac{1}{\beta-1},\frac{1}{\beta-1})^{k}\) when values in \(\R^{k}\) are taken \(\mod\frac{1}{\beta-1}\) component wise.

We note that in the case where \(\beta\) is an integer then \(\beta= 2\) and we are consider a system isomorphic to the usual doubling map on the unit cube. 

\begin{definition}\label{beta rational function}
    A function \(f(\beta)\) is a \(\beta\) rational function if 
    \(f(\beta) \in \Q[\beta]\).
\end{definition}
When the \(\beta\) is clear from context we may simply call the function \(f(\beta)\) or \(f\) a rational function.

We consider maps \(x\to \frac{x+f_{i}}{\beta}\) with \(f_{i} \in \Q[\beta]\).
As \(Q[\beta]\) is a field and \(\Z[\beta]\subset\Q[\beta],\) the following always exists.

\begin{definition}\label{definition: q_{i}}
    Define \(q_{i}\in \Q\) as the minimal element such that such that \(q_{i}f_{i}\in \Z[\beta].\)
\end{definition}

Note that in the case that \(\beta = 2\) the field \(\Q[2]\) is just \(\Q\) and as such rational functions over base \(2\) are just rational numbers. 
Similarly if \(\beta =2\) then for \(f_{i}\in \Q\) then \(q_{i}=\frac{1}{f_{i}}.\)

This allows us to define a key subspace of the \(k+1\) torus. 

\begin{definition} Extended \(\beta\) rational planes \\
    Given an angle vector \(\underline{\theta}= (f_{1},\dots,f_{k})\)with every \(f_{i}\) a \(\beta\) rational function. 
    Define the extended plane at the angle \(\underline{\theta}\) as 
    \[l_{\underline{\theta}} \defeq \{\underline{x} \in [-\frac{1}{\beta-1},\frac{1}{\beta-1})^{k+1}\colon x_{i}\equiv x_{i}' \mod\frac{1}{\beta-1}, \left(\sum_{i=1}^{k}q_{i}x_{i}'\right) - x_{k+1} =0\}.\]
\end{definition}

This is the plane at angle \(\underline{\theta}\) on the \(k+1\) torus. 
Extended planes are invariant under the \(\beta\)ing map on \([-\frac{1}{\beta-1},\frac{1}{\beta-1})^{k}\).
Note that the \(\beta\)ing map is conjugate to the shift map, \(\sigma\), on \(B_{k}^{\N}\) by \(\pi_{\beta}\).

\subsection{Overlapping structure}
The self similar measures with overlaps that we consider have many different and equivalent formulations. 
Two of these have been mentioned in the introduction but  are included here also for completeness.

For \(\beta\) a PV number and \(f_{i}\in \Q[\beta]\) then \(\mu_{\beta,\underline{\theta}}\) is the unique probability measure satisfying,
\[\mu_{\beta,\underline{\theta}}(A) = \frac{1}{k+2}\left(\mu_{\beta,\underline{\theta}}(\beta A) +\mu_{\beta,\underline{\theta}}(\beta A - 1 ) + \sum_{i=1}^{k}\mu_{\beta,\underline{\theta}}(\beta A - f_{i})\right) .\] 

The language of numeration systems gives a good way to refer to these measures as they are the natural measures on the base \(\beta\) numeration systems with the digit set \(D = \{0,1,f_{1},\dots,f_{k}\}\).
As the digit set is a set is it clear we can considering reordering the set without changing the measure or the numeration system associated to the set. 
Formally the measure associated to the \(\beta\) numeration system are the weak star limit measures given in equation \ref{equation : weakstar}.

All of these formulations are clearly related to the IFS formed by the maps \(\{x\to\frac{x+i}{\beta} \colon i \in \{0,1,f_{1},\dots,f_{k}\}\}\). 

A key focus in the study of systems with overlaps are the structural elements known as exact overlaps.

\begin{definition}\label{definition : exact overlaps}
    For an IFS of similarities \(\{g_{1},\dots,g_{k}\}\), if there exists words \(i,j\in \{1,\dots,k\}^{n}\) such that  \(g_{i_{1}}\circ\dots\circ g_{i_{n}} = g_{j_{1}}\circ\dots\circ g_{j_{n}}\) when     
    \(i_{1}\dots i_{n} \neq j_{1}\dots j_{n}\). 
    Then the IFS, maps \(g_{i},g_{j}\) or words \(i,j\) are said to exactly overlaps.

\end{definition}

Exact overlaps are of considerable interest due to the folklore conjecture known as the exact overlap conjecture.

\begin{conjecture}\label{conjecture : exact overlaps}
    For an IFS of similarities of \(\R^{1}\)
    then \(\dim(\mu) = \min\{1,\frac{H(\mathbf{p})}{\lambda({\mathbf{p})}}\}\) unless the IFS has exact overlaps.  
\end{conjecture}

This leads us to introduce the following notions to allow us to count exact overlaps and the growth rate of the number of exact overlaps. 

To count the total number of words in \(A^{n}\) that exactly overlap with a given \(a\in A^{n}\) define 
\[\mathcal{N}_{n}(a,F) \defeq \lvert\{ b \in A^{n} , F_{a} = F_{b}\}\rvert. \]
To count the total number of pairs \(a,b\in A^{n}\) where \(a\) and \(b\) exactly overlap define \[\mathcal{N}_{n}(A,F) \defeq \sum_{a_{1}\dots a_{n}=a\in A^{n}}\mathcal{N}_{n}(a_{1}\dots a_{n},F).\]
This then can be used to find the growth rate of the number of exact overlaps
\[\mathcal{N}(A,F) \defeq \lim_{n\to\infty}\frac{1}{n}\log\mathcal{N}_{n}(A,F).\]

\subsection{Thermodynamic Formalism}

Some notions from Thermodynamic Formalism are used in the conclusions of this work. 
We introduce these notions here and recall the necessary properties when they are used later in the work. 

Before we introduce the potential function we note that there exists a finite number of \(3^{k+1}\times 3^{k+1}\) matrices \(A_{1}\dots A_{3^{k+1}}\) defined independently of choice of \(\beta,f_{1},\dots f_{k}\). 
For details on these matrices see Definition \ref{definition : Markov matrix}.
The following potential function \(\phi\) will later be used to count the growth of exact overlaps.
 \begin{definition}\label{definition : potential intro}
     Let the potential function \(\phi : B_{k+1}^{\N} \to \mathbb{R}\), be given by \[\phi (\underline{a}) = \limsup_{n\to\infty}\log \left( \frac{\mathbf{1}A_{a_{1}\dots a_{n}}(1,0,\dots,0)^{T}}{\mathbf{1}A_{a_{2}\dots a_{n}}(1,0,\dots,0)^{T}}\right).\]
 \end{definition}
This potential function has some key properties.
 \begin{prop}
     There exists a subset \(U \subset B_{k+1}^{\N}\) such that, \\
     \(1)\)  \(U\) has full measure with respect to the fair Bernoulli measure on \(B_{k+1}^{\N}\).\\
     \(2)\)  For any sequence \(\underline{a}\in U\), 
     \[\phi(\underline{a})\vert_{U} = \lim_{n\to\infty}\log\left(\frac{\mathbf{1}A_{a_{1}\dots a_{n}}(1,0,\dots,0)^{T}}{\mathbf{1}A_{a_{2}\dots a_{n}}(1,0,\dots,0)^{T}}\right).\]
 \end{prop}
This statement is proved in a later section, theorem \ref{theorem : well defined}. 

Two key objects in thermodynamic formalism are the potential and the pressure function. 
We are concerned with the pressure of a potential function restricted to extended rational planes in \([-\frac{1}{\beta-1},\frac{1}{\beta-1})^{k+1}\).
Extended rational planes can be given a symbolic representation in the following way.

Define \[X^{n}_{\beta,\underline{\theta}} \defeq \{\underline{x}\in B_{k+1}^{n} \colon \pi_{\beta}(\underline{x}) \in l_{\underline{\theta}}\}.\] 
It is known that \((X^{\N}_{\beta,\underline{\theta}},\sigma)\) is Markov when \(\underline{\theta}\) is rational.
We give a definition of pressure for these subspaces and for a well chosen potential function \(\phi\) defined in \ref{definition : potential}. 
The potential function \(\phi\) is defined on \(B_{k}^{\N}\), we use \(\phi\) as short hand for \(\phi\vert_{X_{\beta,\underline{\theta}}^{\N}}\) to avoid cluttered notation. 

\begin{definition}\label{definition : pressure definition}
        For the space \(X^{\N}_{\beta,\underline{\theta}}\) with the map \(\sigma \colon X^{\N}_{\beta,\underline{\theta}} \to X^{\N}_{\beta,\underline{\theta}}\) and the potential function \(\phi \colon X^{\N}_{\beta,\underline{\theta}} \to \R \) we define the pressure of \(\phi\) on \(X^{\N}_{\beta,\underline{\theta}}\) under \(\sigma\) as 
    \begin{equation*}
        P(X^{\N}_{\beta,\underline{\theta}},\sigma,\phi) \defeq \lim_{n \to \infty}\frac{1}{n} \log\left(\sum_{i_{1}\dots i_{n} \in X^{n}_{\beta,\underline{\theta}}}\exp \left(\sup_{\omega \in [i_{1}\dots i_{n}]}\sum_{j = 0}^{n-1}\phi(\sigma^{j}\omega)\right)\right).
    \end{equation*}
\end{definition}

While it is common practice to write a pressure function only as \(P(\phi)\), we are considering varying spaces so we shall include all parameters for the pressure function. 

\begin{definition}\label{definition : weak gibbs}
A measure \(\nu\) supported on \(X^{\N}_{\beta,\underline{\theta}}\) is called a \emph{Weak Gibbs measure} associated to \(\phi\) if there exists a sequence of positive real numbers \((C_{n})_{n}\) such that \(\lim_{n\to\infty}\frac{\log C_{n}}{n} = 0\) and \(\phi\) such that the following holds, \[\frac{1}{C_{n}} \leq \frac{\mu(\underline{a})}{\exp\left(\sum_{i=0}^{n-1} \left(\phi(\sigma^{i}(a))\right) - n P(X^{\N}_{\beta,\underline{\theta}},\sigma,\phi)\right)} \leq C_{n},\]
for all \( \underline{a} \in X^{\N}_{\beta,\underline{\theta}} \). 
\end{definition}

We now give the definition for a local Weak Gibbs property specialised to our setting. 
For a more general definition and study of local Gibbs properties see \cite{olivier2}.

\begin{definition}\label{definition : local gibbs}
A measure \(\nu\) supported on \([0,\frac{1}{\beta-1})\) is locally Weak Gibbs if the following two statements hold.\\
1) There exists a Weak Gibbs measure \(\eta\) supported on \(E\subset[0,\frac{1}{\beta-1})\) such that \([0,\frac{1}{\beta-1})\setminus E\) is of Hausdorff dimension 0.\\
2)  For any \(y\in [0,\frac{1}{\beta-1}]\setminus E\) 
\[\lim_{r\to0}\left\{\log\nu(B_{r}(x))/\log r\right\} = \alpha \iff \lim_{r\to0}\left\{\log\eta(B_{r}(x))/\log r\right\} = \alpha \]
where \(B_{r}(t)\) denotes the closed ball of radius \(r\) centred at \(t\).
\end{definition}

\section{Motivation of Pressure} \label{section : motivation to pressure}
We begin by motivating the use of a pressure function to study the dimension drop. 
To do this we state the results of Hochman from which we build. 

For an IFS \(G=\{g_{1},\dots g_{l}\}\), let \( r_{i}\) denote the contraction of \(g_{i}\) and \(\mathbf{p}\) a positive probability vector for choosing maps in \(G\). 
Then recall the Lyapunov exponent of \(G\) with probabilities \(\mathbf{p}\) is \(\lambda_{G}(\mathbf{p}) = -\sum_{i}p_{i}\log r_{i} .\)
For all measures \(\nu_{\beta,\underline{\theta}}\) that we are considering the Lyapunov exponent \(\lambda(\mathbf{p})\) is equal to \(\log\beta\). 
We denote the random walk entropy of \(G\) with probabilities \(\mathbf{p}\) as \(h_{RW}(\mathbf{p})\). 
This is currently stated without definition, the definition is given later in Definition~\ref{definition : random walk entropy} when it is used to motivate counting exact overlaps.

Hochman (\!\!\cite{Hochman} Theorem 2.6) states,
\begin{theorem*}\label{theorem : hochman drop}
Let \(G = \{g_{i}\}_{i \in \Lambda}\) be an IFS of similarities in \(\mathbb{R}\) and \(\mathbf{p}\) be a positive probability vector . 
For \(\nu\) the unique probability measure satisfying the self similarity relation for \(G,\mathbf{p}\).
Then \[\dim(\nu) = \min \left\{1,\frac{h_{RW}(\mathbf{p})}{\lambda(\mathbf{p})}\right\},\]
or else \(\min\{d(g_{j},g_{i}) : i,j \in \Lambda^{n}, i\neq j\}\to 0 \) super exponentially.
\end{theorem*}

For the above we note that metric \(d\) isn't in fully generality equivalent to Euclidean distance, which we denote \(\lvert \cdot \rvert \). 
Letting \(g(x) = ax+b, g'(x) = a'x+b'\), Hochman uses the metric \(d(g,g') = \lvert b - b' \rvert + \lvert \log a - \log a' \rvert.\)

Checking that the measures \(\mu_{\beta,\underline{\theta}}\) do not have super-exponential separation has already been done in (\!\!\cite{Hochman_inverse}, Theorem 1.6). 
As the argument is short and used for other purposes later we include it here.
The argument relies on a modification of Garsia separation lemma, \cite{Garsia2}.
\begin{prop}\label{proposition : garsia sep}
    For any PV number \(\beta_{m}\) with conjugates \(\beta_{1},\dots,\beta_{m-1}\)and \(f \in \Q[\beta]\) such that \(f(\beta_{m})\neq 0\) then there exists a constant \(c\) which depends \(\beta,f\) such that 
    \[f(\beta_{i}) > c , \forall i\]
    \begin{proof}
        By definition of a PV number the minimal polynomial for \(\beta\) is a monic polynomial with integer coefficients, \(\lvert\beta_{m}\rvert >1 \) and \(\lvert \beta_{i}\rvert < 1 \) for \(i =1 ,\dots ,m-1.\)
        As \(\Q[\beta] \subset \Z[\beta]\) there exists \(q\in\Z\) such that \(qf \in \Z[\beta]\), then let \(M\) denote the absolute value of the largest coefficient in \(qf\).
        By assumption \(f(\beta_{m})\neq 0\) so we see that \(f(\beta_{i})\neq 0 \forall i\).
        Consequently \(\lvert\prod_{i} qf(\beta_{i})\rvert \geq 1\).
        We obtain \(\lvert qf(\beta_{i})\rvert \leq \frac{M}{1-\lvert\beta_{i}\rvert} \) for \(i = 1,\dots,m-1.\)
        This implies \(qf(\beta_{m}) \geq\prod_{i=1}^{m-1}\frac{1-\lvert\beta_{i}\rvert}{M} \geq 0.\)
        As \(q\) is some fixed finite integer for \(f\) this proves the result.         
    \end{proof}
\end{prop}
This can be extended to any finite collection of \(f_{1},\dots,f_{k}\in \Q[\beta]\) by consider the polynomial given by \(f_{i}-f_{j}\) for all pairs. 

From Garsia's separation lemma, Lemma~\ref{proposition : garsia sep} we can see that there exists a \(c >0\) such that for any \(i,j \in \{1,\dots,l\}^{n}\) then \(d(g_{i},f_{j})\geq \beta^{-n} c.\)
An immediate consequence of this is that the measures \(\mu_{\beta,\underline{\theta}}\) have exponential separation.

Now consider when \(\dim \mu_{\beta,\underline{\theta}}<1\) then by the work of Hochman \(\dim \mu_{\beta,\underline{\theta}} = \frac{h_{RW}(\mathbf{p})}{\log \beta}\). 
Working under the assumption that dimension drop occurs we continue the analysis of \(\frac{h_{RW(\theta)}}{\log \beta}\).
The fact that dimension drop occurs is proven in its own section \ref{section : dim drop}.

Let \(G =\{g_{1},\dots,g_{l}\}\) be an IFS of similarities chosen with positive probabilities \(\mathbf{p}\).
Let \(B_{\mathbf{p}}\) be the Bernoulli measure with probabilities \(\mathbf{p}\). 
The \emph{\(n^{th}\) step random walk entropy} of \(G\) with probabilities \(\mathbf{p}\) is, 
 
\[h_{n}(\mathbf{p}) \defeq -\sum_{a_{1}\dots a_{n} \in \{0,1\}^{n}}B_{\mathbf{p}}([a_{1}\dots a_{n}])\log \sum_{\substack{b_{1}\dots b_{n}\in D^{n} \\ \pi(a_{1}\dots a_{n}) = \pi(b_{1} \dots b_{n})}} B_{\mathbf{p}}\left([b_{1}\dots b_{n}]\right).\]
Then the \emph{random walk entropy} of \(G\) with probabilities \(\mathbf{p}\) is,
\[h_{RW}(\mathbf{p}) \defeq\lim_{n\to\infty}\frac{1}{n}h_{n}(\mathbf{p}).\label{definition : random walk entropy}\]

Throughout this work we are selecting the digits of \(D\) with even probability and so \(\mathbf{p}\) would be the fair \(\frac{1}{k+2}\) vector denoted \(\mathbf{\frac{1}{k+2}}.\)

Akiyama, Feng, Kempton and Persson (\!\!\cite{Akiyama} proposition 3.5) makes a connection between random walk entropy and growth rate of the number of exact overlaps. 
We now modify this argument for our purposes. 
By Jensen's inequality, we obtain 

\begin{align*}
    H_{n}(\mathbf{\frac{1}{k+2}}) & = -\sum_{a_{1}\dots a_{n} \in \D^{n}}B_{(\frac{1}{k+2},\dots,\frac{1}{k+2})}([a_{1}\dots a_{n}])\log \\
    &\qquad\qquad\qquad\qquad\left(\sum_{\substack{b_{1}\dots b_{n}\in D^{n} \\ \pi(a_{1}\dots a_{n}) = \pi(b_{1} \dots b_{n})}}\left( B_{(\frac{1}{k+2},\dots,\frac{1}{k+2})}[b_{1}\dots b_{n}]\right)\right) \\
    & \geq - \log \sum_{a_{1}\dots a_{n} \in \{0,1\}^{n}}(k+2)^{-n} \sum_{\substack{b_{1}\dots b_{n}\in D^{n} \\ \pi(a_{1}\dots a_{n}) = \pi(b_{1} \dots b_{n})}} (k+2)^{-n} \\
    & \geq -\log \sum_{a_{1}\dots a_{n} \in \{0,1\}^{n}}(k+2)^{-2n}\lvert\{b_{1}\dots b_{n} \in D^{n} \colon \pi(a_{1}\dots a_{n})=\pi(b_{1}\dots b_{n}) \}\rvert \\
    & \geq \log(k+2)^{2n} -  \log\mathcal{N}_{n}(D,\{x\to\frac{x+i}{\beta} \colon i \in D\}). 
\end{align*}
This implies that
\begin{equation*}\label{equation : N above n}
    \frac{ h_{RW}(\theta)}{\log\beta}  \geq \frac{2\log k+2}{\log\beta} - \frac{\mathcal{N}(D,\{x\to\frac{x+i}{\beta} \colon i \in D\})}{\log\beta}
\end{equation*}

Under the assumption that dimension drop does occur, we have reduced the the problem of dimension drop to understanding \(\mathcal{N}(D,\{x\to\frac{x+i}{\beta} \colon i \in D\})\). 
We seek to find a dynamical system whose topological pressure can provide and upper bound for \(\mathcal{N}(D,\{x\to\frac{x+i}{\beta} \colon i \in D\})\).
This upper bound would be a lower bound the dimension of the measure.

\section{Construction of Dynamics}\label{section : construction of dynamics}

This section is dedicated to the construction of a dynamical system the topological pressure of which gives and upper bound for the growth rate of the number of exact overlaps in \(\mu_{\beta,\underline{\theta}}\). 
This construction begins with understanding the distance between pairs of words from \(D^{*}\) and the effect that changing this distance has on the number of pairs with that distance. 
We then construct matrices which count the number of pairs of words in \(D^{*}\) with specified distances between members of the pair.
Finally we show that extended rational planes correspond to pairs of points which exactly overlap. 
This allows us to consider only the extended rational planes to count all pair of words that exactly overlap.

To construct this dynamical system we begin by classifying pairs of finite words based upon whether there exist possible extensions which could lead their embedding to be distance zero i.e. exactly overlap. This is done through considering renormalised differences between words.  

\begin{definition}\label{definition : recoverability function}
Let \(\beta \in (1,2]\) a PV number,  $\underline{x},\underline{y} \in D^{\N}$. The \(\beta\) recoverability function is
    \begin{equation*}
        R_{\beta,n}(\underline{x},\underline{y}) = \sum_{i=1}^{n} \beta^{n-i}(x_{i}-y_{i}).
    \end{equation*}
\end{definition}

At times we may suppress the \(\beta\) in \(R_{\beta,n}\), this is done for the ease cumbersome notation when the \(\beta\) is not key to argumentation. 
In the cases when \(\beta\) is suppressed it can be considered an arbitrary and constant PV number.

For two finite words \(x,y\in D^{n}\) should there exist \( x',y'\in D^{m}\) such that \(R_{\beta,n+m}(xx',yy') = 0\) then the pair \(x,y\) could lead to exact overlaps, with \(xx',yy'\).
Note that \(R_{\beta,n+1}(xx_{n+1},yy_{n+1}) = \beta R_{\beta,n}(x,y) +(x_{n+1}-y_{n+1})\) and the maximal digit of \(D\) is \(1\).
Therefore should \(\lvert R_{\beta,n}(x,y)\rvert\geq \frac{1}{\beta-1}\) then \(\lvert R_{\beta,n+1}(xx_{n},yy_{n})\rvert\geq \frac{1}{\beta-1}\) this would mean there cannot exist \(x',y'\) such that \(R_{\beta,n+m}(xx',yy')=0\), by recursive application of the relation.
This gives a bound for the magnitude of \(R_{\beta,n}(x,y)\) such that a possible extension of the pair \(x,y\) exists that leads to an exactly overlapping pair. 

\begin{definition}\label{definition : recoverability}
    For \(x,y\) in \(D^{n}\), call \(x,y\) a \emph{\(\beta\) recoverable pair} if \(\lvert R_{\beta,n}(x,y)\rvert < \frac{1}{\beta-1}.\)
    For \(\underline{x},\underline{y}\) in \(D^{\N}\), call \(\underline{x},\underline{y}\) a \(\beta\) \emph{recoverable pair} if \(\lvert R_{\beta,n}(x_{1}\dots x_{n},y_{1}\dots y_{n}) \rvert < \frac{1}{\beta-1}\) for all \(n\in \N\).\\
    A pair of words or sequences are \(\beta\) \emph{irrecoverable} if they are not recoverable.
\end{definition}

With the aim of understanding the number of recoverable pair of length \(n+1\) in terms of length \(n\) pairs we understand \(R_{n+1}(xx_{n+1},yy_{n+1})\) in terms of \(R_{n}(x,y).\)

\begin{lemma}\label{lemma : r real extentsion lem}
Let \(x,y \in D^{n+1}\times D^{n+1}\) then \(x_{n+1}\) and \(y_{n+1}\) uniquely determine the value of \(R_{\beta,n+1}(x,y)\) in terms of \(R_{\beta,n}(x,y)\). 

\[R_{\beta,n+1}(x,y) =  \begin{cases} 
      \beta R_{\beta,n}(x,y) & x_{n+1} = y_{n+1} \\
      \beta R_{\beta,n}(x,y) - 1 & x_{n+1} = 0 ,\quad y_{n+1} = 1\\
      \beta R_{\beta,n}(x,y) + 1 & x_{n+1} = 1 ,\quad y_{n+1} = 0\\
      \beta R_{\beta,n}(x,y) - f_{i} & x_{n+1} = 0 ,\quad y_{n+1} = f_{i} \\
      \beta R_{\beta,n}(x,y) + f_{i} & x_{n+1}= f_{i} ,\quad y_{n+1} =0 \\
      \beta R_{\beta,n}(x,y) + f_{i} - 1 & x_{n+1} = f_{i} ,\quad y_{n+1} = 1 \\
      \beta R_{\beta,n}(x,y) - f_{i} + 1 & x_{n+1} = 1 ,\quad y_{n+1}=f_{i}\\
      \beta R_{\beta,n}(x,y) +f_{i}-f_{j} & x_{n+1} = f_{i} ,\quad y_{n+1} = f_{j} 
   \end{cases}
\]
\begin{proof}
    This is immediate from a recursive application of the definition of the recoverability function, Definition~\ref{definition : recoverability function}.
\end{proof}
\end{lemma}

Under the assumption that \(R_{\beta,0}(x,y) = 0 \), i.e. the empty word is distance zero from itself, we can construct \(R_{\beta,n}(x,y)\) for any \(n \in \N\) and  any \(x,y \in D^{n} \times D^{n}\) through repeated applications of the above.

The number of exact overlaps can be related to pairs of recoverable words. 
The number of exact overlaps of length \(n\in \N\) for a given \(a \in D^{n}\) can be characterised as  \(\mathcal{N}_{n}(a,\{x\to\frac{x+i}{\beta}
\colon i\in D\}) = \lvert \{ b \in D^{n} , R_{\beta,n}(a,b)=0\} \rvert \). 
From this, it is immediate that for \(a\in D^{n}\), \(\mathcal{N}_{n}(a) \leq \lvert \{ b \in D^{n} , (a,b) \text{ is a \(\beta\) recoverable pair} \}\rvert \).

We now show that recoverable pairs are not just sufficiently close together that there could be an extension leading them to exactly overlap but in fact there always exists such an extension.
This is proven in two part first for considering \(\beta=2\) then considering, \(\beta \in (1,2)\) a PV number. 
The arguments required for each are different as every integer has a finite base \(2\) expansion while not ever integer has finite base \(\beta\) expansion. 

\begin{lemma}\label{lemma : base 2 connection}
    Let \(\beta =2 \) and \(n\geq 1\), \(x,y\in D^{n}\) a recoverable pair. 
    Then there exists \(x',y'\in D^{m} \) for some \(m\geq 0\) such that \(R_{\beta,n+m}(xx',yy')=0.\)
\end{lemma}
    \begin{proof}
        We note that for any recoverable pair \(x,y\in D^{n}\) such that \(R_{2,n}(x,y)=a\) then by Lemma~\ref{lemma : r real extentsion lem} it follows that \(a = (a_{1},\dots,a_{k+1})(f_{1},\dots,f_{k},1)^{T}\).
        We shall assume without loss of generality that \(a >0 \) and \(a_{1},\dots,a_{k}\) are positive. 
        To prove the case that \(a<0 \) the following argument works up to the negation of key signs.
        This shall be proven by combining two claims in turn. 
        Claim 1, there exists \(x',y' \in D^{m'}\) such that \(R_{2,n+m'}(xx',yy')=a' = (0,a'_{2},\dots,a'_{k},a'_{k+1})(f_{1},\dots,f_{k},1)^{T}\).
        Claim 2, for \(s,t\in D^{n'}\) such that \(R_{2,n'}(s,t)=b\) where \(b = (0,\dots,0,b_{k},b_{k+1})\) then there exists \(s',t'\in D^{m}\) such that \(R_{2,n'+m}(ss',tt')=0.\)
        Combining these claims would give that any \((a_{1},\dots,a_{k+1})\) could be reduced to \((0,\dots,0,a_{k},a_{k+1})\) by repeated applications of claim 1, then \((0,\dots,0,a_{k},a_{k+1})\) can be reduced to 0, by claim 2.
        Proof of claim 1:\\
        By definition of \(\Q\) for each \(f_{i}\) then \(\frac{1}{f_{i}}\in \Q\) and \(\frac{1}{f_{i}}f_{i}=1.\) Therefore each \(a_{i}, i\neq k+1\) can be written as a value less than \(\frac{1}{f_{i}}\) by subtracting \(\frac{1}{f_{i}}\) from \(a_{i}\) and adding \(1\) to \(a_{k+1}.\)
        Therefore to reduce \(a'_{i}\) to after appending some \(x',y'\) to \(x,y\) then \(a'_{i} =1/f_{i}\), noting that \(f_{i}\in \Q\) is co-prime then  \(p\) generates the ring of integers mod \(q\) and so any given \(a_{i}\) can be moved to \(1/f_{i}\) as required.      
        Proof of claim 2:\\
        If \(f_{k}=p/q\) is even then \(x'=0,y'=f_{k}\) will reduce \(f_{k}/2\) to zero 
        If \(f_{k}=p/q\) is odd then \(x'=f_{k},y'=1\) will reduce \(f_{k}/2+1/2\) to zero.
        Noting that \(f_{i}\in \Q\) is co-prime then  \(p\) generates the ring of integers mod \(q\) and so any given \(b\) can be moved to \(f_{k}/2\) or \(f_{k}/2+1/2\) as required.
        
    \end{proof}

We now combine this result with the analogous result for \(\beta \in (1,2)\) and show that all recoverable pairs can be extended by some finite pair of words to exactly overlap.

\begin{lemma}\label{lemma : strong connection them}
    Let \(\beta\in(1,2]\) be a PV number and \(n\geq 1\), \(x,y\in D^{n}\) a recoverable pair.
    Then there exists \(x',y'\in D^{m}\) for some \(m\geq 0 \) such that \(R_{\beta,n+m}(xx',yy') =0 .\) 
    \begin{proof}
        Let \(\beta\in (1,2)\) be a PV number then by Garsia's separation lemma, Proposition~\ref{proposition : garsia sep} \(R_{\beta,n}(x,y)\) takes a finite number of values in the range \((-\frac{1}{\beta-1},\frac{1}{\beta-1})\), moreover there exist a \(c>0\) such \(\lvert R_{\beta,n}(x,y) - R_{\beta,n}(x',y')\rvert > c\) whenever \(R_{\beta,n}(x,y) \neq R_{\beta,n}(x',y').\)
        We will show for a given \(R_{\beta,n}(x,y)=a\neq 0\) for some \(x,y\in D^{n}\) that there always exists \(x',y'\in D\) such that \(R_{\beta,n+1}(xx',yy')\) is closed to 0 than \(a\). 
        As \(a\) only takes a finite number of values upon iterating this process, in finite time, will make \(R_{\beta,n+1}(xx',yy')\) closer to \(0\) than the smallest non zero value \(c\) as such \(R_{\beta,n+1}(xx',yy')=0.\)
        We shall be considering only positive \(a\) positive \(R_{\beta,n+1}(xx',yy')=0\) but this argument up to negative of signs holds for negative \(a\).
        For \((R_{\beta,n}(x,y)=a>0\) then \(R_{\beta,n+1}(x',yy')=\beta a + x'-y'\).
        For \(\beta a +x'-y' \) to be closer to 0 than \(a\) then \(\beta a +x'-y' < a\) which implies \(a(\beta-1)< x'-y'.\)
        Noting that \(a \in [c,\frac{1}{\beta-1})\) then such \(x',y'\) can be found. 
        Finally noting that \(\beta a +x'-y' \geq 0\) if and only if \(\beta a \geq y'-x'\) which can be simultaneously satisfied with the above, as \(\beta a \geq y'-x' \iff -\beta a \leq x'-y'\).
        Combing this argument with Lemma~\ref{lemma : base 2 connection} the result is proven.
    \end{proof}
\end{lemma}

We now work to provide a characterisation the dynamics of  \(R_{\beta,n}\to R_{\beta,n+1}\) in terms of a recursive definition that captures how transitioning from \(R_{\beta,n}(x,y)\) to \(R_{\beta,n+1}(x,y)\) effects \(\pi_{\beta}^{-1}(R_{\beta,n}(x,y))\in B_{k+1}^{n}.\)

In the following we consider concatenation of words over \(\{-1,0,1\}\times\dots\times\{-1,0,1\}\) we denote this by writing the letters next to each other without any symbol.

\begin{definition}\label{definition : tilde RN}
     For \(x,y \in D_{k}^{n}\) the function \(\tilde{R}_{\beta,n} : D^{n}\times D^{n} \to (\{-1,0,1\}\times\dots\times\{-1,0,1\})^{n} \) is defined recursively by
\[\tilde{R}_{\beta,n}(x,y) =  \begin{cases} 
      \tilde{R}_{\beta,n-1}(x,y)(\underline{0}) & x_{n} = y_{n} \\ 
      \tilde{R}_{\beta,n-1}(x,y)(0,\dots,0,1)  & x_{n} = 1 ,\quad y_{n} = 0\\
      \tilde{R}_{\beta,n-1}(x,y)(0,\dots,0,-1)& x_{n} = 0 ,\quad y_{n} = 1\\
      \tilde{R}_{\beta,n-1}(x,y)(0,\dots,1,\dots,0)  & x_{n} = f_{i} ,\quad y_{n} = 0 \\
      \tilde{R}_{\beta,n-1}(x,y)(0,\dots,-1,\dots,0) & x_{n}=0 , \quad y_{n}=f_{i} \\
      \tilde{R}_{\beta,n-1}(x,y)(0,\dots,1,0,\dots,-1)  & x_{n} = f_{i} ,\quad y_{n} = 1 \\
      \tilde{R}_{\beta,n-1}(x,y)(0,\dots,-1,0,\dots,1) & x_{n} = 1 ,\quad y_{n} = f_{i} \\
      \tilde{R}_{\beta,n-1}(x,y)(0,\dots,1,0,\dots,-1,0,\dots,0) & x_{n} = f_{i} ,\quad y_{n} = f_{j} ,\quad i\neq j.
   \end{cases}
\]
Finally define \(\tilde{R}_{\beta,0}(x,y) = \underline{0}.\) 
\end{definition}

This \(\tilde{R}_{\beta,n}\) faithfully represents \(R_{\beta}\), as seen in the following lemma. 

\begin{lemma}\label{lemma : r tilde rn equiv }
For every \(n\in \N\) and \(x,y \in D^{n}\) then
\[R_{n}(x,y) = \pi(\tilde{R}_{n}(x,y))(f_{1},\dots,f_{k},1)^{T}.\]
\begin{proof}
    We consider the case that \(n=0\) and continue inductively. 
    If \(n=0\) then we have that \(R_{\beta,0}(x,y)=0\) and that \(\tilde{R}_{\beta,0}(x,y) = \underline{0}, \forall x,y\) by definition. 
    Then \(\pi(\underline{0}) = \underline{0}\) so \(R_{\beta,0}(x,y)=\pi(\tilde{R}_{\beta,0}(x,y))(f_{1},\dots,f_{k},1)^{T}, \forall x,y\).
    
    Assume that \(R_{\beta,n-1}(x,y) = \pi(\tilde{R}_{\beta,n-1}(x,y))(f_{1},\dots,f_{k},1)^{T}\). 
    Then \(R_{\beta,n}(x,y)= \beta R_{\beta,n-1}(x,y)+(x_{n}-y_{n})\) which can be expressed as \( \beta R_{\beta,n-1}(x,y) + R_{\beta,1}(x_{n},y_{n})\).
    Additionally,
    \[\tilde{R}_{\beta,n}(x,y) = \tilde{R}_{\beta,n-1}(x,y)(0,\dots,1,0,\dots,-1,0,\dots,0),\] where the vector \((0,\dots,1,0,\dots,-1,0,\dots,0)\) has \(1\) and \(-1\) in the \(i^{th}\) and \(j^{th}\) component respectively. 
    This expression for \(\tilde{R}_{n}(x,y)\) is also expressible as \(\tilde{R}_{\beta,n}(x,y) = \tilde{R}_{\beta,n}(x,y)\tilde{R}_{\beta,1}(x_{n},y_{n}).\)
    As \[\pi_{\beta}(\tilde{R}_{\beta,n}(x,y)\tilde{R}_{\beta,1}(x_{n},y_{n})) = \beta\pi_{\beta}(\tilde{R}_{\beta,n}(x,y)) + \pi_{\beta}(\tilde{R}_{\beta,1}(x_{n},y_{n})).\] 
    We can apply the assumed identities for lesser \(n\) to see,
    \begin{align*}
        R_{\beta,n}(x,y)&=\beta R_{\beta,n-1}(x,y)+(x_{n}-y_{n}) \\
        &= (\beta\pi_{\beta}(\tilde{R}_{\beta,n}(x,y)) + \pi(\tilde{R}_{\beta,1}(x_{n},y_{n})))(f_{1},\dots,f_{k},1)^{T}\\
        &= \pi_{\beta}(\tilde{R}_{\beta,n}(x,y))(f_{1},\dots,f_{k},1)^{T}. 
    \end{align*}
\end{proof}
\end{lemma}

The possible extensions of \(\tilde{R}_{\beta,n}(x,y)\) give rise to the subclass of words in \(B^{n}_{k+1}\) which can be expressed by as \(\tilde{R}_{\beta,n}(x,y)\) for some \(x,y\in D^{n}\times D^{n}.\) 

\begin{definition}\label{definition : dynamically expressible}
    Any word \(x\in B^{n}_{k+1}\) is called dynamically expressible if the following two properties hold \\
    1) \(\sum_{i=1}^{k+1}\lvert x_{j,i}\rvert \leq 2\)  \\
    2)\(\sum_{i=1}^{k+1}x_{j,i}\leq 1\) for all \(i\leq n \).
\end{definition}

The subset of \(B_{k+1}\) which is all dynamically expressible words will be denoted \(DB_{k+1}\).

This definition is informed by the properties of words expressible by \(\hat{R}_{\beta,n}(x,y)\). 
The first of these properties comes from \(\hat{R}_{\beta,n}\) taking in a pair of words and as such only having two non zero entries in each place.
The second is by virtue of \(\hat{R}_{\beta,n}(x,y)\) appending positive one to the appropriate index according to \(x\) and appending negative one according to \(y\).

\begin{theorem}\label{theorem : all expressible}
    For \(\beta \in (1,2]\) a PV number and every \(n\in \N\). 
    Then every point \(i = R_{\beta,n}(x,y)\) with \(\lvert i\rvert < \frac{1}{\beta-1}\) for some \(x,y \in D^{n}, n\in\N \) is recoverable and dynamically expressible.
    \begin{proof}
        The fact that every \(i = R_{\beta,n}(x,y)\) is dynamically expressible is an immediate consequence of Lemma~\ref{lemma : r tilde rn equiv} and the definition of being dynamically expressible, Definition~\ref{definition : dynamically expressible}. 
        Every value \(\lvert i\rvert < \frac{1}{\beta-1}\) being recoverable is by definition of a pair being recoverable if \(\lvert R_{\beta,n}(x,y)\rvert <\frac{1}{\beta -1}\).
    \end{proof}
\end{theorem}

This establishes that every \(i = R_{\beta,n}(x,y)\) can be expressed as a vector \(z\in B_{k+1}^{l}\) which can be realised by an \(x',y' \in D^{l}\). 

The \(\tilde{R}\) function gives an understanding of how scaled distances can change in relation to k+1 dimensional base \(\beta\) expansions of the distance. 
We use this to construct matrices which count the number of pairs of words with a given distance between them based upon the k+1 dimensional base \(\beta\) expansion of the distance.
As \(B_{k+1}\) is a set it could be re-ordered arbitrarily.
This would give any constructed matrices, a dependency on the ordering of \(B_{k+1}\) for the ordering of its rows and columns. 
Therefore we define \(\mathcal{B} = ((-1,\dots,-1),\dots,(0,\dots,0),(0,\dots,0,1),\dots,(1,0,\dots,0),\dots,(1,\dots,1))\) 
This is the set \(B\) fixed according to the short lex order of its elements.

Before giving the technical details on the construction of the matrices used to define the potential function we seek to provide the intuition behind their construction. 
This is heavily based upon the overlapping structure of cylinder sets in \(B_{k+1}\) when embedded in \(\R^{k+1}\) using the projection map \(\pi_{\beta}\). 

We consider the case of \(k=1\) and the embedding in \(\R^{2}\) for demonstrative purposes but note that the concepts readily generalise to higher dimensional cubes. 
The cylinder set \([0,0]\) is embedded to \([\frac{1}{\beta-1},\frac{1}{\beta-1}]^{2}\).
Assuming that \(\beta \neq 2\) then this cylinder set has overlap with \([1,0]\) which is the set \([1-\frac{1}{\beta-1},1+\frac{1}{\beta-1}]\times[\frac{1}{\beta-1},\frac{1}{\beta-1}]\) in \(\R^{2}\). 
The vertical strip \([1-\frac{1}{\beta-1},\frac{1}{\beta-1}]\times[-\frac{1}{\beta-1},\frac{1}{\beta-1}]\) belongs to both sets, therefore any extensions of \([0,0]\) that fall in this strip will also be representable by some extension of \([1,0]\), which could lead to double counting any possible overlaps. 
Similarly if we were to consider any given corner of the square \([\frac{1}{\beta-1},\frac{1}{\beta-1}]^{2}\) then the corners will overlap with each adjacent square and be counted with multiplicity four. 
This observation can be considered to be that the cylinder \([0,0]\) will overlap with all the cylinders which differ by each coordinate by only one, either plus or minus, i.e the cylinder \([0,0]\) will overlap with \([1,-1]\) but not with \([2,0]\). 

We consider the matrices that are used to define the potential function in two different manners. 
The first is in terms of a matrix of \(1\)'s and \(0\)'s this is to gain an understanding of the topological support for the other form of the matrix. 
The second is a matrix with positive real entries that are used to bound the number of exact overlaps or bound the number of words with a given difference. 
The topological formulation is given first.

\begin{definition}\label{definition : Markov matrix topo}
    For \(\beta\in (1,2]\) a PV number, \(k\geq 1\) and \(a\in B_{k+1}\), define the \(3^{k+1}\times 3^{k+1}\) matrix \(\mathcal{A}_{\beta,a}\) by 
    \[\mathcal{A}_{\beta,a} \defeq \begin{cases}
        1 &\text{if }\tilde{R}_{\beta,n}(x,y) - \tilde{R}_{\beta,n}(s,t) = \mathcal{B}(i) , \tilde{R}_{\beta,n}(x,y)a -\tilde{R}_{\beta,n+1}(ss',tt') \cap [\mathcal{B}(j)]\neq \phi,\\
        0 &\text{otherwise }
    \end{cases}. \]
\end{definition}

We now give the non topological formulation, this is most easily stated in terms of a product of a function that reflects the overlapping structure of cubes in \(\pi_{\beta}(B_{k+1}).\)
This weighting function can be considered to be a version of a Minkowski Steiner formula with where tubular neighbourhoods are taken with respect to the \(1\) metric.
Note that here \(\mathbbm{1}_{E}\) represents the indicator function on the set \(E\), i.e the function \(\mathbbm{1}_{E}(x)= \begin{cases}
    1 \text{ if } x\in E\\
    0 \text{ if } x\notin E.
\end{cases}\) 

\begin{definition}\label{definiton : weighting function}
    Let \(\beta\in (1,2]\) be a PV number, \(k\geq 1\) then \(W\colon B_{k+1}\times B_{k+1}\to \R,\) 
    \[W_{\beta}(x,y) \defeq (1+(k+1)\mathbbm{1}_{x=y}) \prod_{i=1}^{k+1}\left( (\mathbbm{1}_{\lvert x_{i}-y_{i}\rvert=0}\frac{2}{\beta-1})+(\mathbbm{1}_{\lvert x_{i}-y_{i}\rvert=1}\frac{2}{\beta-1}-1)\right)\]
\end{definition}

We note that the division by two in the above weighting can be seen to account for what would other wise be double counting caused by the symmetry of \(W\) on a pair.
The indicator function multiplying by \(k+1\) is to represent that there is more ways to write \(0\) in \(B_{k+1}-B_{k+1}\) namely the diagonal in this set.

This weighting function can then be used to construct a matrix that counts the number of words with a given distance between them.
The construction uses the idea that the weighting function represents the fraction of all possible extensions that overlap a given position so weighting that by the number of ways to write such an extensions gives the desired count. 

\begin{definition}\label{definition : markov matrix}
    For \(\beta \in (1,2]\) a PV number, \(k\geq 1\), \(a\in B_{k+1}\), define the \(3^{k+1}\times3^{k+1}\) matrix \(A_{\beta,a}\) by
    \[A_{\beta,a}(i,j)\defeq \sum_{x\in DB_{k+1}}W_{\beta}(ix,\alpha+j).\]
\end{definition}


These matrices have two distinct positive entry structure based on if \(\beta\in (1,2)\) a PV number or if \(\beta = 2.\)
This structural difference arises from the inherent overlap between the cylinders \([0],[1]\) for PV numbers and no such overlap in base 2. 
This means that the weighting function is not required for base \(2\) and the matrix \(A_{\beta,a}\) can be stated simply in terms of \(\tilde{R}_{2,n}\).

\begin{lemma}\label{lemma : markov matrix base 2}
        For \(k\geq 1\) and \(a \in B_{k+1}\), let 
        \begin{align*}
        E= \{s',t'\in D\times D\colon s,t \in D^{n},& \tilde{R}_{2,n}(x,y) - \tilde{R}_{2,n}(s,t) = \mathcal{B}(i) , \\
        &\tilde{R}_{2,n}(x,y)a -\tilde{R}_{2,n+1}(ss',tt') =\mathcal{B}(j) \},
        \end{align*}
        then \[A_{2,a}(i,j) = \lvert E \rvert .\]
    \begin{proof}
        By Definition~\ref{definiton : weighting function} we see that \(W_{2}(x,y) = (1+(k+1)\mathbbm{1}_{x=y}) \prod_{i=1}^{k+1}\left( \mathbbm{1}_{\lvert x_{i}-y_{i}\rvert=0}\right)\)
        Now consider \(\sum_{x\in DB_{k+1}}W_{2}(ix,\alpha j)\) as in Definition~\ref{definition : markov matrix}, we see that 
        \begin{align*}
            \sum_{x\in DB_{k+1}}W_{2}(ix,\alpha j) &= 
            \sum_{x\in DB_{k+1}}W_{2}(ix,\alpha j) \\
            &= (1+(k+1)\mathbbm{1}_{ix=\alpha+j}) \prod_{i=1}^{k+1}\left( \mathbbm{1}_{\lvert ix_{i}-(\alpha+j)_{i}\rvert=0}\right)\\
            & =\begin{cases}
                k+2 &\text{ if } ix=\alpha+ j\\
                1 &\text{ if} 
            \end{cases}
        \end{align*}
        Note that up expressing as cardinality of a set as opposed to a sum the claim is proven.
    \end{proof}
\end{lemma}

With the definition of the \(A_{\beta,a}\) matrices simplified in the case of \(\beta=2\) we now consider \(\beta \in (1,2)\) a PV number.
In this case it can be shown that the matrices \(A_{\beta,a}\) have the same positive entry structure for all such \(\beta\) and each \(a\). 

\begin{lemma}\label{lemma : markov matrix positive structure}
Let \(\beta,\beta'\in (1,2)\) be a PV numbers such that \(\beta \neq \beta'\),\(k\geq 1\) and \(a\in B_{k+1}\) then \(A_{\beta,a}(i,j) \neq 0 \iff A_{\beta',a}(i,j) \neq 0\). 
    \begin{proof}
        This is immediate from noting that \(\frac{2}{\beta-1}\) and \(\frac{2}{\beta-1}-1\) are both strictly positive for all such \(\beta\) as such \(W_{\beta}(ix,\alpha j)\) will non zero if and only if \(W_{\beta}(ix,\alpha j)\) is non zero.
    \end{proof}
\end{lemma}

This motivates the suppression of \(\beta\) in \( A_{\beta,a}\) unless \(\beta =2\).

We define a Markov process on the elements of \(B_{k+1}\).
While this process is well defined, the overlaps intrinsic to PV numbers mean that multiple elements of \(B_{k+1}^{n}\) will embed to the same point in \([-\frac{1}{\beta-1},\frac{1}{\beta-1})^{k+1}\). 
This Markov process will give an upper bound for the number of pairs of words with a given distance between them based upon the \(k+1\) dimensional base \(\beta\) expansion of the distance.

\begin{lemma}\label{lemma : count equiv}
    Let \(\beta\in (1,2]\) a PV number, \(k\geq1\) \(f_{1},\dots,f_{k}\in \Q[\beta]\) and \(x,y \in D^{n}\) such that \(\tilde{R}_{\beta,n}(x,y) = z\) then,
    \[\lvert \{ x,y \in D^{n} \colon  R_{\beta,n}(x,y)= \pi_{\beta}(z)(f_{1},\dots,f_{k},1)^{T}\}\rvert \leq \mathbf{1}A_{\underline{z}}(1,0,\dots,0)^{T}  .\]
    \begin{proof}
        The Markov process is a stationary and as such it suffices to show that the claim holds for a matrix product of length \(1\).
        By Lemma~\ref{lemma : markov matrix base 2} gives that the matrix \(A_{z}\) is equivalent to counting extensions with the prescribed difference. 
        The the vectors \((1,\dots,1), (1,0,\dots,0)\) sum over all values with difference \((0,\dots,0)\) after application of \(A_{z}\).
        
        Now we consider the case \(\beta\) a PV number.
        The matrices \(A_{\beta,a}\) will count the same quantity as \(W_{\beta}(ix,\alpha j)\) summed over all dynamically expressible words in \(B_{k+1}\).

        Recalling that the definition of dynamically expressible is inspired by the possible values \(\tilde{R}_{\beta}(s,t)\) can take, we see that summing over all dynamically expressible words is equivalent to considering all possible extensions for a pair of words in \(B_{k+1}\times B_{k+1}\).
        Therefore it is sufficient to show that \(W_{\beta}(ix,\alpha j)\) gives the number of \(x \in B_{k+1}\) such that \(\tilde{R}_{\beta}(s,t) = i\) and \(\tilde{R}_{\beta}(sx,t\alpha) = j.\)
        Noting that \(\tilde{R}_{\beta}(sx,t\alpha)= \tilde{R}_{\beta}(s,t)\tilde{R}_{\beta}(x,a) = i\tilde{R}_{\beta}(x,\alpha) = i(x-\alpha)\) then \(i(x-\alpha) =j\) so \(ix = \alpha + j\). 
        By construction \(W_{\beta}(ix,\alpha+j)\) counts the number of \(x\) such that the claim follows.
        
        
    \end{proof}    
\end{lemma}

This allows us to understand exact overlaps through the matrices \(A_{z}.\)

\begin{corollary}\label{corollary : exact overlap symb}
    For \(\beta \in (1,2]\) a PV number, \(k\geq 1 \), \(f_{1},\dots,f_{k}\in \Q[\beta]\), \(n\in \N\) and \(Z = \{ z \in B_{k+1}^{n} \colon \pi_{\beta}({z})(f_{1},\dots,f_{k},1)^{T} = 0\}\).
    Then 
    \[ \sum_{z\in Z} \mathbf{1}A_{z}(1,0,\dots,0)^{T} \geq\sum_{ a \in D^{n}} \mathcal{N}_{n}(a) = \mathcal{N}_{n}\left(D,\left\{x\to\frac{x+i}{\beta} \colon i \in D\right\}\right). \]
    \begin{proof}
        This is an immediate consequence of \(\mathcal{N}_{n}(a) = \lvert \{b \in D^{n}, R_{\beta,n}(a,b) = 0 \} \rvert\) and Lemma~\ref{lemma : count equiv}.
    \end{proof}
\end{corollary}

Sub-spaces of \([-\frac{1}{\beta-1},\frac{1}{\beta-1})^{k+1}\) which are invariant under the dynamics of the \(\beta\)ing map, and correspond to the set of pairs of sequences which contain exact overlaps are key to the main results of this work. 

\begin{definition}\label{definition : binary cube}
    For \(\beta\in (1,2]\) a PV number,\(k\geq 1\) then a \emph{binary cube} is a subset of \([-\frac{1}{\beta-1},\frac{1}{k+1})^{k+1}.\)
    A binary cube has side length \(2\beta^{-n}\) and  centre at \((\frac{i_{1}}{\beta^{n}},\dots,\frac{i_{k+1}}{\beta^{n}})\), where \[\pi_{\beta}^{-1}\left(\frac{i_{1}}{\beta^{n}},\dots,\frac{i_{k+1}}{\beta^{n}}\right) \in B_{k+1}^{n}.\]
\end{definition}

We give \([-\frac{1}{\beta-1},\frac{1}{\beta-1})^{k+1}\) the usual dynamics of the \(\beta\)ing map. 
Given a binary cube \(A\), we see that \(T_{\beta}^{-1}(A)\) is again a binary cube or \([-\frac{1}{\beta-1},\frac{1}{\beta-1})^{k+1}\).

Binary cubes can be seen as the geometric analogue of the symbolic cylinder sets. 
We now define the geometric version of pairs corresponding to exact overlaps of length \(n\). 
By Corollary \ref{corollary : exact overlap symb} this would be all points \(x \in [-\frac{1}{\beta-1},\frac{1}{\beta-1})^{k+1}\) such that \(x(f_{1},\dots,f_{k},-1)^{T}=0\). 

Recall the definition of an extended \(\beta\) rational plane.
For \(f_{i} \in \Q[\beta]\) and \(q_{i}\in \Q\) where \(q_{i}\) is the least such value that \(q_{i}f_{i}\in \Z[\beta].\)
Then given an angle vector \(\underline{\theta}= (f_{1},\dots,f_{k})\) define the extended plane at the angle \(\underline{\theta}\) as 
\[l_{\underline{\theta}} \defeq  \left\{\underline{x} \in [-\frac{1}{\beta-1},\frac{1}{\beta-1})^{k+1}\colon x_{i}\equiv x_{i}' \mod\frac{1}{\beta-1}, \left(\sum_{i=1}^{k}q_{i}x_{i}'\right) - x_{k+1} =0\right\}.\]

From this definition a few properties are immediate. 
\begin{itemize}
    \item \(l_{\underline{\theta}}\) is invariant under \(T_{\beta}\), \(T^{-1}(l_{\underline{\theta}}) = l_{\underline{\theta}}\)
    \item If \(x \in l_{\underline{\theta}}\) then there exists a \(n \in \N\) such that \[T_{\beta}^{n}(x) = x' \text{ where } x'(f_{1},\dots,f_{k},-1)^{T}= 0.\]
\end{itemize}

\section{Potential and Pressure}\label{section: potential and pressure}

\subsection{Potential}
With the link between exact overlaps and sub-spaces of \([-\frac{1}{\beta-1},\frac{1}{\beta-1})^{k+1}\) established, we now turn our attention to a potential function on \([-\frac{1}{\beta-1},\frac{1}{\beta-1})^{k+1}\). 

The pressure of this potential function on extended \(\beta\) rational planes captures the maximal growth rate of the number of exact overlaps. 
This section is motivated by the works of Chazottes and Ugalde \cite{chazottes} and begins by introducing the basic objects of their work. 
In this, we are working with more general sequences and cannot use their techniques to recover any Gibb's properties.  

\begin{definition}\label{definition : potential}
For a sequence \(\underline{z} \in B_{k+1}^{\N}\) define \(\phi \colon B_{k+1}^{\N} \to \mathbb{R}\) by 
    \[\phi(\underline{z}) \defeq \limsup_{n\to \infty} \log\left(\frac{\mathbf{1}A_{z_{1}}A_{z_{2}}...A_{z_{n}}(1,0,\dots,0)^{T}}{\mathbf{1}A_{z_{2}}A_{z_{3}}...A_{z_{n}}(1,0,\dots,0)^{T}}\right).\]
\end{definition}

To understand \(\phi\), we introduce open simplexes upon which the matrices, \(A_{i} , i \in B_{k+1}\), can act. This allows us to gain an understanding of the variation of \(\phi\).

\begin{definition}\label{definition : open simplex}
    Let \(n \geq 1\) then \(E_{n}\) denote the open \(n\) simplex, also called the \(n\) simplex for ease, defined by
    \begin{equation*}
        E_{n} \defeq \left\{ \underline{x} \in \mathbb{R}^{n} \colon (x_{1},\dots,x_{n})\in(0,1)^{n} , \sum_{i=1}^{n=1}x_i = 1 \right\}. 
    \end{equation*}
    Moreover for \(I = \{i_{1},\dots,i_{l}\}\) let \(E_{n-l,I}\) denote the sub simplex formed by removing the coordinates in \(I\),
    \begin{equation*}
        E_{n-l,I} \defeq \left\{ \underline{x} \in \mathbb{R}^{n} \colon I=\{1,\dots,n\}, i\in I\setminus\{i_{1},\dots,i_{l}\},  x_{i} \in(0,1) ,x_{i}=0 ,  \sum_{i=1}^{n}x_i = 1 .\right\} 
    \end{equation*}
\end{definition}

The simplexes \(E_{n-1,i}\) are the open faces of the \(n\) simplex \(E_{n}\) and as such \(E_{n-1,i} \not\subset E_{k}.\) 
Note that, \(\overline{E_{k}}\) the closure of \(E_{k}\) allows the \(x_{i} = 0 \text{ or }1\), similarly for \(E_{k-l,\{i_{1},\dots,i_{l}\}}\) .

An open \(n\) simplex can be acted on by an appropriately normalised \(n\times n\) matrix. 
Properties of such an action under the Hilbert metric allows us to conclude where \(\phi\) is well defined, 

\begin{definition}\label{definition : normalised matrix action}
For \(n\geq 1\) and a \(n\times n\) non negative matrix \(M\), the normalised matrix action \(F_{M}(\underline{x}):\overline{E_{n}}\to\overline{E_{n}}\) is
    \begin{gather}
        F_{M}(\underline{x}) \defeq \frac{M\underline{x}}{\lVert M\underline{x} \rVert}.
    \end{gather}
\end{definition}

The open \(n\) simplex is now given the Hilbert metric, this is chosen for its resemblance to the potential function \(\psi\).

\begin{definition}\label{definition : hilbert metric}
For \(n\geq 1\), \(x,y \in E_{n}\), the Hilbert metric is,
    \begin{equation*}
        d_{E_{n}}(\underline{x},\underline{y}) \defeq \log \left( \frac{\max_{1\leq i \leq n} \frac{x_{i}}{y_{i}}}{\min_{1 \leq i \leq n}\frac{x_{i}}{y_{i}}} \right).
    \end{equation*}
For \(I = \{i_{1},\dots,i_{l}\}, J= \{j_{1},\dots,j_{l}\}\) then\(x \in E_{n-l,I},y\in E_{n-l,J}\) then, 
    \begin{equation*}
        d_{E_{n-l,I,J}}(\underline{x},\underline{y}) \defeq \log\left( \frac{\max\{\frac{x_{a}}{y_{a}},\frac{x_{j_{b}}}{y_{i_{b}}}\colon a \in \{1,\dots,n\}\setminus(I\cup J),\ 1\leq b\leq l\}}{\min\{\frac{x_{a}}{y_{a}},\frac{x_{j_{b}}}{y_{i_{b}}}\colon a \in \{1,\dots,n\}\setminus(I\cup J),\ 1\leq b\leq l\}}\right).
    \end{equation*}
\end{definition}

The Hilbert metrics allow us to understand the distance between points on simplexes of the same dimension. 
This means that we can consider the distance between \(x,y \in \{z\in[0,1]^{n}\colon\  \sum_{i=1}^{n}z_{i}=1\} \), i.e \(x,y\) which have the same number of zero entries. 
This is as the Hilbert metric treats \(\delta E_{n}\) as a boundary at infinity and we have noted \(E_{n-l,I} \subset \delta E_{n}\).
For an in depth justification of the Hilbert metric and background on this topic see \cite{Kohlberg}.

The potential function \(\psi\) is defined in terms of products of \(A_{i}\)'s. 
To understand a general \(n \times n\) matrix product \(M_{\underline{i}}= M_{1}\dots M_{k}\) we consider the iterative effect of each \(M_{\underline{i}}\) on \(E_{n}\).
Therefore we must ensure that the correct Hilbert metric is being used at each stage, this is done through placing restrictions on matrices we consider. 
This restriction is found later, Theorem~\ref{theorem : contractive matirx}. 
To this end we introduce the following notion, which can be though of as a weakening of primitivity. 

We say the \(j^{th}\) row of a matrix, \(M\), is positive if \(M(i,j) > 0\) for all \(i\).
Moreover, we say \(M\) has \(k\) positive rows if \(M\) has \(k\)  rows which are positive.
Similarly, we say the \(j^{th}\) row of a matrix, \(M\), is a zero row if \(M(i,j) = 0\) for all \(i\).
We say \(M\) has \(k\) zero rows if \(M\) has \(k\) rows that are zero rows. 
Let \(\pos(M)\) denote the number of positive rows of a matrix \(M\) and \(\zero(M)\) denote the number of zero rows of a matrix \(M\).
In the case that \(\zero(M)>0\) let \(\overline{\zero(M)}\) denote the set of indexes for the zero rows.

For a \(n\times n\) matrix \(M\), \(x,y\in E_{n}\) let, 
\[d_{E_{F_{M}}}(F_{M}(x),F_{M}(y)) \defeq d_{E_{pos(M),\overline{\zero(M)} ,\overline{\zero(M)} }}(F_{M}(x),F_{M}(y))
.\]
In the case that \(\pos(M)=n\) we see in this case simply recovers the usual Hilbert metric on the n simplex, \(d_{E_{n}}.\)

Given an understanding of which Hilbert metric for apply to a given matrix we can move to understanding the contractivity of a matrix. 
Denote the matrix products \( A_{z} = A_{1}\dots A_{n} \) for \(z \in B_{k+1}^{n}\) and \(A_{\underline{z}} = A_{1}A_{2}\dots \) for \(\underline{z}\in B_{k+1}^{\N}\). 

\begin{definition}\label{definition : contractive matrix}
    For \(n\geq 1\) and a \(n\times n\) matrix, \(M\) we call \(M\) contractive if, 
    
    \[\sup_{\underline{x},\underline{y} \in E_{n}} \left( \frac{d_{E_{F_{M}}}(F_{M}(\underline{x}),F_{M}(\underline{y}))}{d_{E_{n}}(\underline{x},\underline{y})} \right) < 1.\]
    
    A word \(z \in B_{k+1}^{n}\) is called contractive when the matrices \(A_{z}\) is contractive. \\
    A contractive sequence \(\underline{z}\in B_{k+1}^{\N}\) is called infinitely contractive if it contracts \(E_{3^{k+1}}\), the open \(3^{k+1}\) simplex, to a single point. 
\end{definition}

We now seek to express contractive and infinitely contractive in terms of properties of the matrices \(A_{i}, i \in B_{k+1}\). 
This is done so that we can find an uniform contraction coefficient for contractive matrices. 

We state a theorem of Chazottes and Ugalde \cite{chazottes} in the terms of this paper. 

\begin{theorem}\label{theorem : contractive matirx}
    Let \(n\geq 1\) and \(M\) be a non-negative \(n\times n\) matrix with \(1<j \leq n\) positive rows and \(n-j\) zero rows. 
    Then A is contractive.
    \begin{proof}
        This is the statement of (\!\!\cite{chazottes} Lemma 2) where we consider the positive rectangular matrices as the positive rectangular matrices defined by removing the \(n-j\) rows of zeros. 
    \end{proof}
\end{theorem}

The contractivity of \(M\) in the above theorem in part relies up the metric \(d_{n,I,J}(M(x),M(y))\) being able to be applied for all \(x,y \in E_{n}\). 
As the matrix \(M\) is positive on all rows it is not zero we see that \(M(x)\in E_{\pos(M)}\) for all \(x\in E_{n}\) and so it is well posed to ask \(d_{N-l,I,J}(M(x),M(y))\) for all \(x,y\in E_{n}\). 

We now specialise to consider nine specific matrices and their action on \(E_{9}\).
This is done to establish the measure of the set of infinitely contractive matrices. 
We note these are the matrices the matrices given by Definition~\ref{definition : markov matrix} for \(k=1\).

By Theorem~\ref{theorem : contractive matirx}, contractivity of a matrix is understood in terms of being strictly positive on any rows which are not zero rows.
The following results categorise which products of \(A_{i}, i\in B_{2}\) have such a structure.

\begin{lemma}\label{lemma : max 1 zero row.}
    Let \(\beta \in (1,2]\) be a PV number then any finite product of the matrices $A_{i}, i \in B_{2}$ is non zero. 
    \begin{proof}
        This is proven through a large case wise argument considering the relative positions of positive entries of matrices being taken in the product.
        Fix \(\beta\) throughout the proof, any \(\beta\) in subscripts will be suppressed for readability and because all argumentation works for any given \(\beta.\)
        We do not present the entire argument here but a subset of the cases to show the arguments structure.
        Begin by considering products of the form \({A}_{(0,0)}{A}_{(1,x)}\) for \(x\in\{-1,0,1\}\). 
        As \({A}_{(0,0)}\) is positive for \({A}_{(0,0)}(1,1),{A}_{(0,0)}(2,2)\) and \(A_{(1,x)}(2,2)\) is positive for all \(x\) the product will have positive entry in position \((2,2)\).
        Continuing in such a fashion for all remaining positive structures in \({A}_{a}\) and any products found completes the proof. 
        This process is finite as every product of length three has at least one positive row. 
    \end{proof}
\end{lemma}

A consequence of this non degeneracy condition that the number of contractive matrices of a given length can be understood.

\begin{lemma}\label{lemma : contractivity dens rough}
    For \(\beta\in (1,2]\) a PV number, \(k= 1\), \(f_{1}\in \Q[\beta]\) then at least \(78\) of the \(729 = 9^{3}\) of the matrices in \(\{A_{\beta,a_{1}}A_{\beta,a_{2}}A_{\beta,a_{3}}\colon a_{1},a_{2},a_{3}\in B_{2}\}\) are contractive.
    \begin{proof}
    This follows from an exhaustive calculation of products of length 3 and the application of Theorem \ref{theorem : contractive matirx} to the length 3 products.
    \end{proof}
\end{lemma}

This yields a finite number of contractive matrices of length three which we use to define a universal contraction coefficient. 
This coefficient is the lower bound for the amount of contraction that occurs under the action of a contractive matrix.
We consider matrices which are not contractive to be like isometries of \(E_{9}\). 
This is due to the fact that the matrices will preserve distances between some pair of points and so will have a contraction coefficient of 1.

\begin{definition}\label{definition : contraction coefficient}
    For words \(\underline{\omega} \in  (B_{2})^{3}\) the maximal contraction coefficient is \[\tau =  \max\{c : c\neq 1 , c\sup_{\underline{x},\underline{y}\in E_{9}}d_{E_{9}}(\underline{x},\underline{y}) \geq \sup_{\underline{x},\underline{y}\in E_{9}}d_{E_{A_{\omega}}}(A_{\omega}\underline{x}
,A_{\omega}\underline{y}) \}.\]
\end{definition}

We now consider the size of the set of infinitely contractive matrices \(A_{i},i\in B_{2}^{\N}.\)
We show that this is a full measure subset of \(B_{2}^{\N}\) and later use pull backs of this set to establish properties of \(\phi.\)

\begin{theorem}\label{theorem : full mes set}
    The set of infinitely contractive sequences, 
    \(\mathcal{I}\subset B_{2}^{\N}\) is full measure with respect to the Bernoulli(\(1/9,\dots,1/9\)) measure on \(B_{2}^{\N}\).
    \begin{proof}
        Consider re-coding \(B_{2}^{\N}\) by elements of \(B_{2}^{3}\). 
        By Lemma~\ref{lemma : contractivity dens rough}, we know that at least 78 elements of \(B_{2}^{3}\) are contractive. 
        Therefore under the re-coding by \((B_{2}^{3})^{\N}\), at least \(\frac{78}{729}\) of the elements are contractive. 
        Re-coding \((B_{2}^{3})^{\N}\) again by \(\{c,i\}\) according to whether the element of \(B_{2}^{3}\) is contractive or not, we obtain \(\{c,i\}^{\N}\) with the Bernoulli \((\frac{78}{729},\frac{651}{729})\) measure.
        Applying the strong law of large numbers to the set \(\{\underline{\omega}^{\N},         \lim_{n\to\infty}\frac{1}{n}(\sum_{i} \mathbbm{1}_{c}\omega_{i}) = \frac{78}{729}\}\), we see that \(\mu(\{\underline{\omega}^{\N}, \lim_{n\to\infty}\frac{1}{n}(\sum_{i} \mathbbm{1}_{c}\omega_{i}) = \frac{78}{729}\}) = 1\). 
        The measure of finitely contractive sequences is 0 as an immediate consequence.
    \end{proof}
\end{theorem}

With the special case of \(\mathcal{I}\) established for \(B_{2}\) now show how the matrices \(A_{i},i\in B_{2}\) occur as substructures in \(A_{i},i\in B_{k+1}.\) 
We then use this to show that \(\phi\) is well defined on a full measure subset of \(B^{\N}_{k+1}\). 

To do this we introduce notation for the minor of a matrix.
For \(a \in B_{k+1}\) let \(A_{\beta,a}\vert_{i,j}\) be the \(9\times9\) matrix taken by considering only the rows and columns of \(A_{\beta,a}\) for which \((\mathcal{B}(l)_{i},\mathcal{B}(l)_{j})\in B_{2}\) and \(\mathcal{B}(l)_{m}=0\) for \(m\neq i, j\) and all \(\mathcal{B}\). 

\begin{lemma}\label{lemma: a9 substrucutre}
    For \(\beta \in (1,2]\) a PV number, \(k\geq 1 \) and \(a\in B_{k+1}\) with \(i,j \in \{1,\dots,k+1\}\)  and \(a' \in B_{2}\) such that \(a'_{1} = a_{i}, a'_{2}=a_{j}\) then, 
    \[A_{\beta,a}\vert_{i,j}(x,y)\neq0 \iff A_{\beta,a'}(x,y)\neq 0.\]
    \begin{proof}    
        By definition of \(A_{\beta,a}\), Definition~\ref{definition : markov matrix} 
        \(A_{\beta,a}\vert_{i,j}(x,y)\neq0 \iff A_{\beta,a'}(x,y)\neq 0\) is equivalent to \( \sum_{x\in DB_{k+1}}W_{\beta}(ix,\alpha+j) \neq 0 \iff  \sum_{x\in DB_{k+1}}W_{\beta'}(ix,\alpha+j)\).
        As the \(\frac{2}{\beta-1},\frac{2}{\beta-1}-1\) are both non zero for all \(\beta \in (1,2)\) then the non zero nature of both sums is entirely dependant on \(\lvert x_{i}-y_{i}\rvert\). 
        Noting that \(\lvert x_{i}-y_{i}\rvert\) is independent of \(\beta\) the claim follows. 
    
    \end{proof}
\end{lemma}
We now use the presence of these substructures and the full measure of \(\mathcal{I}\) with respect to \(B_{2}^{\N}\), to show that \(\phi\) is defined on a full measure subset of \(B_{k+1}^{\N}\) for arbitrary \(k.\) 

\begin{theorem}\label{theorem : well defined}
    For \(k\geq 1\), the function \(\phi\) is well defined on a full measure set of \([-\frac{1}{\beta-1},\frac{1}{\beta-1})^{k+1}.\)
    \begin{proof}
        First we note some properties of the set \(\mathcal{I}\). 
        We can represent \(\mathcal{I}\) as \(\mathcal{I} =(\mathcal{I}_{l},\mathcal{I}_{r})\) is symmetric so for \((x,y)\in \mathcal{I} \implies (y,x)\in \mathcal{I}\). 
        So we can find a \(\mathcal{I}_{s}=\mathcal{I}_{l}=\mathcal{I}_{r}\) such that \(\mathcal{I}= \mathcal{I}_{s} \times \mathcal{I}_{s} \setminus C\) for some \(C\) which is measure 0. 
        We will write \(C = (C_{l},C_{r})\). 

        The function \(\phi\) will be well defined if the associated matrix product for each point contracts the entire simplex to a point. 
        To contract the entire simplex to a single point we note that it must contract every 
        sub space that is isomorphic to the open nine simplex to a  point.
        Therefore the set of associated points where \(\phi\) is well defined must project to give \(\mathcal{I}\) on any subspace isomorphic to a nine simplex.
        
        Now we define \(P_{i,j}(\hat{\mathcal{I}}) = (0,\dots,0,a_{i},0,\dots,0,a_{j},0,\dots,0)\). 
        So \\ \(P^{-1}_{i,j} = \{(\underline{\alpha_{1}},\dots,\underline{\alpha_{n}})\in B_{n}^{\N}, (\underline{\alpha_{i}},\underline{\alpha_{j}})\in \mathcal{I}\}.\)
        Now consider 
        \begin{equation*}
        \hat{\mathcal{I}} = \bigcap_{\substack{i \\ j \neq i}}P_{i,j}^{-1}\mathcal{I} .
        \end{equation*}
        \begin{align*}
            P^{-1}_{i,j}(\mathcal{I}) &\bigcap P_{i,k}^{-1}(\mathcal{I}) = \{(\underline{\alpha_{1}},\dots,\underline{\alpha_{n}})\in B_{n}^{\N}, (\underline{\alpha_{i}},\underline{\alpha_{j}})\in \mathcal{I}, (\underline{\alpha_{i}},\underline{\alpha_{k}})\in \mathcal{I}\} \\
            &= \{(\underline{\alpha_{1}},\dots,\underline{\alpha_{n}}), p\in\{1,\dots,n\}\setminus\{i,j,k\}, \underline{\alpha_{p}}\in B_{1},  \\
            & \quad \quad(\underline{\alpha_{i}},\underline{\alpha_{j}},\underline{\alpha_{k}})\in \mathcal{I}_{s}\times \mathcal{I}_{s} \times \mathcal{I}_{s} \setminus C_{l}\times C_{r} \times B_{1} \setminus C_{l} \times B_{1} \times C_{r}\}
        \end{align*}
        Then 
        \begin{align*}
        \bigcap_{j}P^{-1}_{i,j}(A) = 
        &\mathcal{I}_{s}\times\dots\times \mathcal{I}_{s} \setminus C_{l}\times C_{r} \times B_{1} \times \dots \times B_{1} \\
        &\setminus C_{l}\times B_{1}\times C_{r} \times B_{1} \times \dots \times B_{1} \setminus \dots \setminus c_{l}\times B_{1} \times \dots B_{1} \times C_{r}  
        \end{align*}

        This is the finite removal of measure zero sets and so \(\hat{\mathcal{I}}\) is full measure.
    \end{proof}
\end{theorem}

From this we are able to conclude when \(\psi\) is well defined. 

\begin{corollary}\label{corollary : lim sup equiv}
    \(\phi(\underline{x})\) is well defined on a full measure set, namely, all \(x\in B^{\N}_{k+1}\) that are infinity contractive.  
    \begin{proof}
        Define \(\psi(v)\colon E_{3^{k+1}} \to \R\) by \(\psi(v) \defeq \log(1,1,1,1)A_{x_{0}}v\). 
        The function \(\psi\) is continuos in \(v\in E_{3^{k+1}}\). 
        As the matrices \(A_{i}\) are contractions in \(E_{3^{k+1}}\) space and \(E_{3^{k+1}}\) is complete with respect to the Hilbert metric we can apply Banach fixed-point theorem to \(\psi.\) 
        By definition \(\underline{x}\in B_{k+1}^{\N}\) that is infinitely contractive, map \(E_{3^{k+1}}\) to a point. 
        Therefore \(\underline{x}\in B_{k+1}^{\N}\) that is infinitely contractive have \(\phi(\underline{x})\) as well defined. 
    \end{proof}
\end{corollary}

As \(\phi\) is well defined on a full measure sub set of the \(k+1\) torus we can establish the following relation to growth rate of exact overlaps. 

For the sake of readability let \(D' \defeq \{x\to\frac{x+i}{\beta} \colon i\in D\}\).
\begin{lemma}\label{lemma : sup n equiv}
    Let \(Z_{n} = \{z \in B_{k+1}^{n} \colon \pi({z})(f_{1},\dots,f_{k},-1)^{T} = 0\}\) then, 
    \[\log\frac{\mathcal{N}_{n}(D,D')}{\mathcal{N}_{n-1}(D,D')} \leq  \sum_{x\in Z_{n}}\sup_{y\in [x]}\phi(y) .\]
    \begin{proof}
        Let \(x \in \{ x \in B^{n} :  \pi_{\beta}{(x})(p,-q)^{T} = 0\} \) then by Corollary \ref{corollary : exact overlap symb}, 
        \[ \sum_{x} (1,1,1,1)A_{x}(1,0,0,0)^{T} = \mathcal{N}_{n}(D,D') .\] 
        Dividing through by the same expression for \(n-1\), 
        \begin{align*}
            \frac{\mathcal{N}_{n}(D,D')}{\mathcal{N}_{n-1}(D,D')}
            &= \frac{\sum_{x} (1,1,1,1)A_{x_{1}}\dots A_{x_{n}}(1,0,0,0)^{T}}{\sum_{z} (1,1,1,1)A_{z_{2}}\dots A_{z_{n}}(1,0,0,0)^{T}} \\ 
            &\leq \sum_{x} \frac{(1,1,1,1)A_{x_{1}}\dots A_{x_{n}}(1,0,0,0)^{T}}{(1,1,1,1)A_{x_{2}}\dots A_{x_{n}}(1,0,0,0)^{T}} \\
            &\implies \log \frac{\sum_{a \in D^{n}} \mathcal{N}_{n}(a) }{\sum_{b\in D^{n-1}} \mathcal{N}_{n}(b)}\\
            &\leq  \log \sum_{x} \frac{(1,1,1,1)A_{x_{1}}\dots A_{x_{n}}(1,0,0,0)^{T}}{(1,1,1,1)A_{x_{2}}\dots A_{x_{n}}(1,0,0,0)^{T}} \\
            &\leq  \sum_{x} \log \frac{(1,1,1,1)A_{x_{1}}\dots A_{x_{n}}(1,0,0,0)^{T}}{(1,1,1,1)A_{x_{2}}\dots A_{x_{n}}(1,0,0,0)^{T}} 
            \leq \sum_{x}\sup_{y\in [x]}\phi(y). 
        \end{align*}
    \end{proof}
\end{lemma}

\section{Proof of key theorems}\label{restults sec}

\subsection{Dimension Drop}\label{section : dim drop}
This section is dedicated to showing that this generalisation of the Bernoulli convolution we are considering do exhibit dimension drop.

Due to the exact overlap conjecture, Conjecture~\ref{conj:exact overlaps} we understand that dimension drop for self similar measures occurs in the presence of exact overlaps.
As discussed in the introduction it has been shown that Bernoulli convolutions have exact overlaps for \(\beta\in (1,2)\) a PV number. 
The measures that we are considering can be seen as Bernoulli convolutions with additional digits, \(f_{i}\in \Q[\beta].\)
If every \(f_{i}\) is expressible with a finite code over \(\{0,1\}\) then it is clear that each \(f_{i}\) only induces more exact overlaps, this would lead one to expect that \(\mu_{\beta,\underline{\theta}}\) would still exhibit dimension drop.

It is worth noting that the bounds for dimension drop given in Hochman \cite{Hochman} depend on the number of digits and past a certain number of digits it becomes unclear whether a system should exhibit dimension drop. 

The simplest example of a self similar measure with overlaps is the projection of the Sierpi\'nski triangle onto the measure \(\nu_{2,\{0,1,f_{1}\}}\).
This projected measure exhibits dimension drop if and only if \(f_{1}\in \Q\), this can be seen as an example in \cite{Hochman_inverse}.

This is in part due to the fact that if \(f_{1}\in\R\setminus\Q\) then \(f_{1}\) would not exhibit exact overlaps and the system of \(\{x\to\frac{x+i}{2}, i \in \{0,1,f_{1}\}\}\) would not have exponential separation.

More generally we note that the additional digits having the finite coding is key in this dimension drop else the exponential separation which is key to the understanding of the fractal at all scales breaks. 
This is as the elements of \(\Q[\beta]\) are analogous to the behaviour of elements of \(\Q\) in base \(2\). 

We are unaware of any concise statement of the dimension drop for this generalisation of Bernoulli convolutions in the literature and as such include one here.

\begin{theorem*}
    For \(\beta \in(1,2]\) a PV number, \(k\geq 1\), \(f_{1},\dots,f_{k} \in \Q[\beta]\) and \(\underline{\theta} = \{f_{1},\dots,f_{k}\}\). 
    Then 
    \[\dim(\mu_{\beta,\underline{\theta}}) < 1.\]
\end{theorem*}

We note that if \(\beta=2\) this result is known to be true and can be seen as the consequence of results in a few different areas, for example weak separation theorems or projections of integer lattices. 

This can be proven in two different methods; deriving the result from more general theory, or constructing the result from the properties of this system specifically. 

To prove this from general results we note the following two theorems in the literature.

We state Proposition~3.2 from Py{\"o}r{\"a}l{\"a} \cite{pyor} in the language of this paper. 
\begin{theorem}
    Let \(G = \{g_{1},\dots,g_{k}\}\) be an homogenous IFS of affine contractions on \(\R\).
    If there exists \(c,d > 0 \) such that for \(i, j \in \{1,\dots,k\}^{n}\) then \(\lvert g_{i} - g_{j}\rvert = 0 \) or \(\lvert g_{i} - g_{j}\rvert \geq c d^{-n}\), then for \(\nu\) the unique probability measure associated to \(G\),  \[\dim(\nu)=1.\]
\end{theorem}

While Br{\'e}mont Theorem~2.3~\cite{brémont2020} states the following in the language of this work. 
\begin{theorem}
    Let \(\beta\in (1,2]\) a PV number with \(k\geq 1\) and \(\underline{\theta} = (f_{1},\dots,f_{k})\) with \(f_{i}\in \Q[\beta]\) for all \(i\). 
    Then \(\mu_{\beta,\underline{\theta}}\) is singular with respect to the Lebesgue measure.
\end{theorem}

These two theorems readily combine to give that for this class of self similar measures with overlaps that singularity is equivalent to dimension drop and that all measures in this class are singular. 

The proof of this result using less general results is significantly longer but does have the advantage of constructing a specific automaton in the process.
This automaton reveals interesting and useful structural elements which allow for easier comparison of this work to known literature. 

To prove dimension drop occurs for this class of self similar measures with overlaps without these more general results, we build a series of lemmas which we then combine with classic techniques to show the singularity of a measure.  

For the case that \(\beta = 2\) these properties are known but for the sake of completeness we include a short proof of the result using the same techniques as that for \(\beta\) a PV number. 

Recall Lemma \ref{lemma : zero one} that if \(\beta \in (1,2)\) a PV number then \(\beta\) satisfies a polynomial with coefficients in \(\{-1,0,1\}\).

\subsubsection{Singularity}
We use the techniques of Fourier analysis to establish the singularity of the measures \(\mu_{\beta,\underline{\theta}}\). 
This technique was used by Erd{\"o}s in \cite{Erdos} to show that the Bernoulli convolution is singular with respect to Lebesgue. 

\begin{definition}\label{definition: fourier transform}
For the measure \(\mu_{\beta,\underline{\theta}}\) define the Fourier transform \(\hat{\mu}_{\beta,\underline{\theta}}\) as 
\[\hat{\mu}(\eta)_{\beta,\underline{\theta}}  \defeq \int_{\R}e^{it\eta}d\mu_{\beta,\underline{\theta}}(t).\]
\end{definition}

This allows us to calculate the Fourier transform of the generalised Bernoulli convolutions we consider. 
\begin{prop}
    The Fourier transform of \(\mu_{\beta,\underline{\theta}}\) is 
    \[\hat{\mu}_{\underline{\theta}}(\eta) = \prod_{n=0}^{\infty} \left(\frac{2}{k+2}\cos(\beta^{n}\eta)+\frac{1}{k+2}\sum_{i=1}^{k}\cos((2f_{i}-1)\beta^{n}\eta)\right).\]
    \begin{proof}
        This follows from a standard calculation.
    \end{proof}
\end{prop}

A consequence of the Riemann Lebesgue lemma is that if \(\hat{\mu}_{\beta,\underline{\theta}}(\eta)\not\to 0\) as \(\eta \to \infty \) implies that \(\mu_{\beta,\underline{\theta}}\) is not absolutely continuos with respect to the Lebesgue measure.

The family of measures \(\mu_{\beta,\underline{\theta}}\) as a family of self similar measures are known to be either purely singular or absolutely continuous \cite{Feng}.
Therefore showing lack of absolute continuity suffices, to show this we follow argumentation outlined by Solomyak in \cite{mandelbrot2004fractal}. 

\begin{lemma}\label{lemma : singular}
    For \(\beta \in (1,2]\) a PV number and \(k \geq 1\) with \(f_{1},\dots,f_{k} \in \Q[\beta]\) then \(\mu_{\beta,\underline{\theta}}\) is singular with respect to the Lebesgue measure. 
    \begin{proof}
        Assume \(\beta \in (1,2)\) a PV number. 
        let \(\beta_{2},\dots, \beta_{m(\beta)}\) be the Galois conjugates of \(\beta\), recall that \(\beta_{2},\dots, \beta_{m(\beta)}\) all lie inside the unit disc in \(\mathbb{C}.\)
        Therefore \(\max_{i\geq2}\lvert\beta_{i} \rvert = D \in (0,1).\)
        The functions \(\beta^{n}+\sum_{i=2}^{m(\beta)}\beta_{i}^{n}\) is a symmetric function of the roots of a inter coefficient minimum polynomial. 
        From this symmetry and the algebraic properties of \(\beta\), the function \(\beta^{n}+\sum_{i=2}^{m(\beta)}\beta_{i}^{n} = \alpha \in \N\) for every \(n\in\N\). 
       
        Now consider the Hausdorff distance, \(\dist\) between \(\beta^{n}\) and \(\Z\). \\
        Using \(\beta^{n}+\sum_{i=2}^{m(\beta)}\beta_{i}^{n} = \alpha \in \N\) and \(\max_{i\geq2}\lvert\beta_{i} \rvert = D \in (0,1)\) then, \[\dist(\beta^{n},\Z) \leq (m(\beta)-1)D^{n}.\]\label{hausdorff distance eq}

        With this distance property established we now consider a chosen sequence of \(\eta_{n}\) such that \(\hat{\mu}_{\beta,\underline{\theta}}\not\to 0\) as \(\eta_{n}\to\infty\).

        As \(\beta^{n}\) is an algebraic integer and \(\beta \neq 2\)  then \(\beta^{n}\not\in \frac{1}{2}\Z\).
        This implies that \(\hat{\mu}_{\beta,\underline{\theta}}(\pi \beta^{N})\neq 0\).
        
        For \(N>1\),
        \begin{align*}
            \hat{\mu}_{\beta,\underline{\theta}}(\beta^{N}\pi) = &\left(\frac{2}{k+2}\cos(\beta^{N}\pi)+\frac{1}{k+2}\sum_{i=1}^{k}\cos((2f_{i}-1)\beta^{N}\pi) \right) \\
            &\qquad\qquad\left(\frac{2}{k+2}\cos(\beta^{N-1}\pi)+\frac{1}{k+2}\sum_{i=1}^{k}\cos((2f_{i}-1)\beta^{N-1}\pi) \right)\\
            &\dots \left(\frac{2}{k+2}\cos(\pi)+\frac{1}{k+2}\sum_{i=1}^{k}\cos((2f_{i}-1)\pi) \right) \\
             &\qquad\qquad\left(\frac{2}{k+2}\cos(\beta^{-1}\pi)+\frac{1}{k+2}\sum_{i=1}^{k}\cos((2f_{i}-1)\beta^{-1} \pi) \right)   \dots \\
            &= \prod_{n=1}^{N}\left(\frac{2}{k+2}\cos(\beta^{n}\pi)+\frac{1}{k+2}\sum_{i=1}^{k}\cos(\beta^{n}(2f_{i}-1)\pi) \right)\ \hat{\mu}_{\beta,\underline{\theta}}(\pi).
        \end{align*}
            From this we begin to consider the magnitude of the Fourier transform along this sequence. 
            In the following we use the Hausdorff distance inequality established in equation \ref{hausdorff distance eq}
        \begin{align*}
            \lvert \hat{\mu}_{\beta,\underline{\theta}}(\beta^{N}\pi)\rvert =& \left\lvert \prod_{n=1}^{N}\left(\frac{2}{k+2}\cos(\beta^{n}\pi)+\frac{1}{k+2}\sum_{i=1}^{k}\cos(\beta^{n}(2f_{i}-1)\pi) \right)\ \hat{\mu}_{\beta,\underline{\theta}}(\pi)\right\rvert \\ 
            &\geq \left\lvert \prod_{n=1}^{\infty}\left(\frac{2}{k+2}\cos(\beta^{n}\pi)+\frac{1}{k+2}\sum_{i=1}^{k}\cos(\beta^{n}(2f_{i}-1)\pi) \right)\ \right\rvert\left\lvert \hat{\mu}_{\beta,\underline{\theta}}(\pi)\right\rvert \\
            & \geq \biggl\lvert \prod_{n=1}^{\infty}\biggl(\frac{2}{k+2}\cos((m(\beta)-1)D^{n}\pi)+ \\
            &\qquad\qquad\quad \frac{1}{k+2}\sum_{i=1}^{k}\cos((m(\beta)-1)D^{n}(2f_{i}-1)\pi) \biggr)\ \biggr\rvert\left\lvert \hat{\mu}_{\beta,\underline{\theta}}(\pi)\right\rvert .
            \
        \end{align*}
        
        As the \(f_{i}'s\) are a finite collection of no zero elements of \(\Q[\beta]\) then there exists a \(f_{m}\) such that \(\lvert 2f_{m}-1 \rvert \) is minimal. 
        For the sake of notational ease let \(f_{m}' = 2f_{m}-1\).
        Therefore,

        \begin{align*}
            \lvert \hat{\mu}_{\beta,\underline{\theta}}(\beta^{N}\pi)\rvert \geq& \biggl\lvert \prod_{n=1}^{\infty}\biggl(\frac{2}{k+2}\cos((m(\beta)-1)D^{n}f_{m}'\pi)+ \\ 
            & \qquad\qquad\qquad\frac{1}{k+2}\sum_{i=1}^{k}\cos((m(\beta)-1)D^{n}f_{m}'\pi) \biggr)\ \biggr\rvert\left\lvert \hat{\mu}_{\beta,\underline{\theta}}(\pi)\right\rvert .
        \end{align*}
        
        Finally for this case note that 
        \begin{align*}
        \biggl\lvert \prod_{n=1}^{\infty}\biggl(\frac{2}{k+2}\cos(&(m(\beta)-1) D^{n}f_{m}'\pi)+ \\
        &\frac{1}{k+2}\sum_{i=1}^{k}\cos((m(\beta)-1)D^{n}f_{m}'\pi) \biggr)\ \biggr\rvert\left\lvert \hat{\mu}_{\beta,\underline{\theta}}(\pi)\right\rvert \not\to 0 .
        \end{align*}
        As \(D^{n}\to 0 \) as \(n\to\infty\).

        In the case that \(\beta = 2\) we see that every \(2f_{i}-1 \in \Q\) so denote \(f'_{i} = (2f_{i}-1)\). 
        Each \(f_{i}'\) will have a denominator, as \(2f_{i}-1 \in \Q\), let \(q_{i}\) be this denominator. 
        The Fourier transform for this case is \(\hat{\mu}_{2,\underline{\theta}} = \left(\frac{2}{k+2}\cos(2^{n}\eta)+\frac{1}{k+2}\sum_{i=1}^{k}\cos(f_{i}'2^{n}\eta)\right)\). 
        Taking the sequence \(\eta_{n}= \prod_{i} q_{i}\pi\) we see that \(\hat{\mu}_{2,\underline{\theta}} = 1 \) for all n, so \(\hat{\mu}_{2,\underline{\theta}}\to 1 \neq 0\) as required. 
    \end{proof}
\end{lemma}

This shows that for all generalised Bernoulli convolutions we consider are singular. 
It now remains to show that these measures are dynamically invariant to be able to show dimension drop occurs, this is done in the next section. 

\subsubsection{Dynamic invariance}\label{section: dynamic invariance}
With the singularity of \(\mu_{\beta,\underline{\theta}}\) established we need to now show that \(\mu_{\beta,\underline{\theta}}\) is equivalent to a dynamically invariant measure. 
We follow the ideas introduced in Vershik and Sidorov \cite{Vershik} to construct a dynamically invariant measure on an automaton. 

This begins by defining an automaton.
\begin{definition}\label{definition : automoton}
    For \(\beta\in(1,2]\) a PV number, \(k\geq 1\), \(\underline{\theta} = (f_{1},\dots,f_{k})\), with \(f_{i}\in\Q[\beta]\) for all \(i\), we define the finite state automaton \((G,E)_{\underline{\theta}}\) as follows 
    \begin{itemize}
        \item Vertex set \(G = \{g \in \Z[\beta,f_{1},\dots,f_{k}], \lvert g\rvert < \frac{1}{\beta-1} \}\)
        \item Edge set \(E = \{ (a,b)\in G\times G, \exists  x,y\in D ,  \beta a + (x-y) = b\} \) 
        \item Label the edge from \(a\) to \(b\) by \((a,b)\) if \((a,b) \in E\)
    \end{itemize} 
\end{definition}

Note that by Garcia separation Proposition~\ref{proposition : garsia sep}, G is finite. 

An automaton is called \emph{strongly connected} if every state is connected to every other state by some path in the automaton. 

\begin{lemma}\label{lemma : automota is strongly conncected}
    The automaton \((G,E)_{\underline{\theta}}\) is strongly connected.
    \begin{proof}
        Note that the set \(G\) is the same as the values that \(R_{\beta,n}(x,y)\) takes in the range \((\frac{-1}{\beta-1},\frac{1}{\beta-1})\). 
        By comparing Definitions~\ref{definition : automoton},~\ref{definition : markov matrix},
        note that there is an edge between \(a,b\) if and only if for \(x,y\in D^{n}\) such that \(R_{\beta,n}(x,y)= a \) there exists \(x',y'\in D\) such that \(R_{\beta,n+1}(xx',yy')=b.\)  
        As every state of \(R_{\beta,n}(x,y)\) is dynamically expressible and that there exists extensions \(x',y'\) such that \(R_{\beta,n'}(xx',yy')=0\) by Lemma~\ref{lemma : betary rep} we can deduce that 0 is connected to every state and every state is connected to 0. 
        Therefore \(G\) is strongly connected.
    \end{proof}
\end{lemma}

We can associate weights to the edges in \(E\). 
Weight the edge from a to b according to the number of pairs \(x,y \in D\times D\) such that \(\beta a + (x-y) = b\). 
Denote such an edge \((a,b)\) and its weight \(\lvert (a,b) \rvert \). 
Call this weighted edge set \(E'\). 
It is clear that \((a,b)\in E \implies (a,b) \in E'\).
Should a pair \((a,b) \in G \times G\) and \((a,b)\notin\ E\) then we give the edge weight 0.

Define the transition matrix \(M^{(G,E')_{\underline{\theta}}}\). 
This matrix gives the probability of transitioning from the state \(i\) to \(j\) in \((G,E')\). 
It is often useful to refer to a specific transition probability or entry of \(M_{(i,j)}^{(G,E')_{\underline{\theta}}}\). 
Given this we write \(\underline{p}_{i,j} = M^{(G,E')_{\underline{\theta}}}_{(i,j)}\).

This formulation of the dynamics of \(R_{\beta,n}(x,y)\) on \((G,E')\) give another equivalent definition of \(\mu_{\beta,\underline{\theta}}.\)

\begin{lemma}\label{lemma : lim def mu theta}
    For \(\beta \in (1,2)\) a PV number with \(k\geq 1\) with \(\underline{\theta} = (f_{1},\dots,f_{k})\) and 
    \(f_{i}\in \Q[\beta]\) for all \(i\) then, 
    \[\mu_{\beta,\underline{\theta}}(A) = \lim_{n\to\infty}\sum_{\substack{a \\ \pi(a) \in A}}\left(\underline{p}_{0,a_{1}}\prod_{i=1}^{n}\underline{p}_{a_{i},a_{i+1}}\right).\]
    \begin{proof}
        For a given \(a \in A\),  \(\underline{p}_{0,a_{1}}\prod_{i=1}^{n}\underline{p}_{a_{i},a_{i+1}}\) gives the probability of all \(b \in D^{n}\) such that \(\lvert a-b\rvert \leq \beta^{-n}\). 
        As \(n \to \infty\) this gives probability of all \(b\) such that \(a=b\). 
        Summing over all \(a\in A\) gives \(\mu_{\beta,\underline{\theta}}(A)\). 
    \end{proof}
\end{lemma}

We note that the measure \(\mu_{\beta,\underline{\theta}}\) is not shift invariant because the unique start vertex forces \(\underline{p}_{0,a_{1}}\) as the first term in Equation~\ref{equation : mu weak star}. 
Crucially, this is the only position at which a choice is forced and so it is the only location where shift invariance fails.

From the automaton \((G,E')_{\underline{\theta}}\) we construct the invariant \(\tilde{\mu}_{\beta,\underline{\theta}}\) as the following weak star limit  
\[\tilde{\mu}_{\beta,\underline{\theta}}(A) \defeq \lim_{n\to\infty} \sum_{\substack{a \\ \pi(a) \in A}}\left( \sum_{\substack{j\in G\\ (j,a_{1} \in E')}}\underline{p}_{j,a_{1}} \left(\prod_{i=1}^{n}\underline{p}_{a_{i},a_{i+1}}\right)\right).\]

This constructed measure \(\tilde{\mu}_{\beta,\underline{\theta}}\) can be readily shown to be dynamically invariant and equivalent to \(\mu_{\beta,\underline{\theta}}.\)
\begin{lemma}\label{lemma: invariant equiv}
    For \(\beta \in (1,2]\) a PV number, \(k\geq 1\), \(\underline{\theta} = (f_{1},\dots,f_{k})\)  with \(f_{i}\in \Q[\beta]\) for all \(i\) then \(\mu_{\beta,\underline{\theta}}\) and \(\tilde{\mu}_{\beta,\underline{\theta}}\) are equivalent as measures.
    \begin{proof}
        By the definition of \(\mu_{\beta,\underline{\theta}}\) and \(\tilde{\mu}_{\beta,\underline{\theta}}\) their formulations differ in only their initial positions, \(\sum_{\substack{j\in G\\ (j,a_{1} \in E')}}\underline{p}_{j,a_{1}}\) and \(\underline{p}_{0,a_{1}}\).
        This unique position of difference lets us show that there exists constants \(c_{l},c_{r}\) such that \(c_{l}\tilde{\mu}_{\beta,\underline{\theta}} \leq \mu_{\beta,\underline{\theta}} \leq c_{r}\tilde{\mu}_{\beta,\underline{\theta}}\). 
        As \(\underline{p}_{(0,a_{1})}\) is a summand of \(\sum_{j}\underline{p}_{(j,a_{1})}\), \(c_{r} = 1 \) is a valid constant. 
        Similarly, if we take \(c_{l} = \min_{\substack{a \\ (0,a)\in E'}}\{\frac{\underline{p}_{(0,a)}}{\sum_{j}\underline{p}_{(j,a)}}\}\) then, \\
        \(c_{l} \tilde{\mu}_{\beta,\underline{\theta}} \leq \mu_{\beta,\underline{\theta}} \leq c_{r}\tilde{\mu}_{\beta,\underline{\theta}}\). 
        All that remains is to show that \(c_{l}\) is non zero. 
        As we only take minimum over vertices connected to \(0\), \(\underline{p}_{(0,a)}\) is positive therefore \(c_{l}\) is non zero.  
    \end{proof}
\end{lemma}

This has now shown that \(\mu_{\beta,\underline{\theta}}\) is equivalent to a dynamically invariant measure. 

We deviate from completing the proof of dimension drop here to see how these automota sit within the literature. 
We compare this construction with the construction in \cite{laulq}, which allows for the calculation of \(L^{2}\) dimension of measures such as \(\mu_{\theta}\).

In \cite{laulq} they calculate they present a method for calculating the \(L^{2}\) of a wide range of self similar measures in \(\R^{d}\). 
The measures we are considering are a subclass of those for which their result holds. 
A key part of their construction is to consider the transition probabilities, as informed by a large matrix, between ``near by'' compositions of functions. 

In \cite{laulq} they show they are able to consider a finite matrix, this is in direct parallel to the finite image of recoverable \(R_{\beta,n}\). 
Similarly they consider taking the irreducible component of their matrix, this is analogous the strongly connected property of the automaton, \((G,E)_{\underline{\theta}}.\) 

From these analogies we can see this method could be seen to recover a specification of \cite{laulq} or give different method to construct the matrixes in this specific case.
As seen in section \ref{section : construction of dynamics} instead of pursue this direct line of inquiry we choose to generalise this construction in such a manner that looses the need to completely re construct the matrix or automaton from scratch for every member of this family.
This is done to give an alternative approach to construct a global lower bound for the dimension of this family of measures. 

\subsubsection{Dimension drop}
We now combine the singularity of \(\mu_{\beta,\theta}\) and equivalence to an invariant measure, \(\tilde{\mu}_{\beta,\underline{\theta}}\), to show that the dimension of \(\mu_{\beta,\theta}\) must be less than 1. 

\begin{theorem}\label{dimdropoc}
     Let \(\beta \in (1,2]\) a PV number ,\(k\geq 1\), \(\underline{\theta} = (f_{1},\dots,f_{k})\) where \(f_{i}\in \Q[\beta]\) for all \(i.\)
     Then \(\dim(\mu_{\beta,\underline{\theta}}) < 1.\)
    \begin{proof}
        It is know that the parry measure is the unique measure of maximal entropy for this class of generalised Bernoulli convolutions as well as a much larger class of symbolic systems. 
        A work that gives this result is \cite{hofbauer1978shifts} which showed that there exists a unique measure of maximal entropy  measure \(\nu\), for \(\beta\in(1,2)\). 
        We have shown in Lemma~\ref{lemma : singular} that \(\mu_{\beta,\underline{\theta}}\) is singular with respect to Lebesgue and in Lemma~\ref{lemma: invariant equiv} that \(\mu_{\beta,\underline{\theta}}\) is invariant with respect to the \(\beta\)ing map. 
        Combining these results by the uniqueness of \(\nu\) and singularity of \(\mu_{\beta,\underline{\theta}}\) then \(\dim(\mu_{\beta,\underline{\theta}}) < 1.\)
    \end{proof}
\end{theorem}

\subsection{Gibbs properties}\label{gibbsprop}

To establish the Gibbs properties of \(\mu_{\beta,\underline{\theta}}\), we follow the ideas of Olivier et al, (\!\!\cite{olivier},\cite{olivier2}) and introduce a related measure \(\bar{\mu}_{\beta,\underline{\theta}}\) the Gibbs structure of which is related to \(\mu_{\beta,\underline{\theta}}\) at a local level. 
For a complete background on local Gibbs properties, the general construction of \(\bar{\mu}_{\beta,\underline{\theta}}\) and its properties see (\!\!\cite{olivier},\cite{olivier2}) and references within. 

We begin by recalling the definition of Weak Gibbs and Locally Weak Gibbs.

\begin{definition*}\label{definition : weak gibbs again}
A measure \(\nu\) supported on \(X\) is called a Weak Gibbs measure associated to \(\phi\) if there exists a sequence of positive real numbers \((C_{n})_{n}\) such that \(\lim_{n\to\infty}\frac{\log C_{n}}{n} = 0\) and \(\phi\) such that the following holds, \[\frac{1}{C_{n}} \leq \frac{\mu(\underline{a})}{\exp\left(\sum_{i=0}^{n-1} \left(\phi(\sigma^{i}(a))\right) - n P(X^{\N}_{\underline{\theta}},\sigma,\phi)\right)} \leq C_{n},\]
for all \( \underline{a} \in X \). 
\end{definition*}

\begin{definition*}\label{definition : localdef2}
A measure \(\nu\) supported on \([0,\frac{1}{\beta-1})\) is locally Weak Gibbs if the following two statements hold.\\
1) There exists a Weak Gibbs measure \(\eta\) supported on \(E\subset[0,\frac{1}{\beta-1})\) such that \([0,\frac{1}{\beta-1})\setminus E\) is of Hausdorff dimension 0.\\
2)For any \(y\in [0,\frac{1}{\beta-1}]\setminus E\) 
\[\lim_{r\to0}\left\{\log\nu(B_{r}(x))/\log r\right\} = \alpha \iff \lim_{r\to0}\left\{\log\eta(B_{r}(x))/\log r\right\} = \alpha,\]
where \(B_{r}(t)\) denotes the closed ball of radius \(r\) centred at \(t\).
\end{definition*}

To define \(\bar{\mu}_{\beta,\underline{\theta}}\) we first introduce additional transformations and an extended digit set to support \(\bar{\mu}_{\beta,\underline{\theta}}\).

For \(i\in D\) define \(F_{\beta,i}\) by \(F_{\beta,i}(x) \defeq \frac{x+i}{\beta}\) so \(F^{-1}_{i}(x) \defeq \beta x - i\). 
Using the maps \(F_{\beta,i}\), and the self similar relation, equation \ref{equation : mu self sim} we see that, 
\[\mu_{\beta,\underline{\theta}}(A + x) = \sum_{j\in D} \frac{1}{k+2} \mu_{\beta,\underline{\theta}}(F^{-1}_{\beta,i}(A)+ \beta x + (i-j)).\] 
Looking at \(\mu_{\beta,\underline{\theta}}(F^{-1}_{\beta,i}(B)+y)\), we see \(\mu_{\beta,\underline{\theta}}(F^{-1}_{\beta,i}(B)+y) = 0\) when \(y \not\in (-1,  1)\). 
When \(y \in (-1,1)\), we write \(x \triangleright y\) if \(\beta x-y + i \in D\) for some \(i \in D\). \\
This allows us to define the expanded digit set \(\mathcal{D}\). 
The advantage of this set is that it allows us to express \(\pi_{\beta}(D_{k}^{n})\) entirely with independent digits of length one, for certain \(n\in\N\). 
Analogously, \(\mathcal{D}\) is the set of values that the recoverability function \(R_{\beta,n}(x,y)\) takes in the interval \((\frac{-1}{\beta-1},\frac{1}{\beta-1})\).
Formally, 
\[\mathcal{D} \defeq \bigcup_{n} \bigcup_{k=0}^{k_{n}} \{y = \beta i_{k} + (i-j) \colon (i,j) \in \{0,1\} \times D, -1 < y < 1\},\]and \(i_{0} =0 \). 
We note a difference in the set \(\mathcal{D}\) and the analogous \(\mathcal{I}_{\beta,d}\) in
\cite{olivier}. 
As we are considering rational digits our set \(\mathcal{D}\) contains more points than (proposition 2.5 \cite{olivier}) states, but is still finite. 

Akin to the construction of transition matrices in \((G,E)\) we can construct the matrices \(M_{0},M_{1}\) which define, in terms of a points base \(\beta\) digit expansion, the transition probabilities inside \(\mathcal{D}\).

For \( i\in \{0,1\}, j,k \in \{0, \dots ,\lvert \mathcal{D}\rvert \} ,i_{k},i_{j} \in \mathcal{D}\) we have
\[M_{i}(h,k) \defeq \lvert\{ (x,y) \in D\times D, x-y = j \colon i + \beta i_{k} - i_{k} = j  \in D\}\rvert\]
Therefore we can see the following relation,
\[\begin{pmatrix}
    \mu_{\beta,\underline{\theta}}(B+i_{0}) \\
    \vdots \\
    \mu_{\beta,\underline{\theta}}(B+i_{k})  
\end{pmatrix} = M_{i}\begin{pmatrix}
    \mu_{\beta,\underline{\theta}}(F_{\beta,i}^{-1}(B)+i_{0}) \\
    \vdots \\
    \mu_{\beta,\underline{\theta}}(F_{\beta,i}^{-1}(B)+i_{k})
\end{pmatrix}.\]

Define the measure \(\bar{\mu}_{\beta,\underline{\theta}}(B) \defeq \frac{\sum_{j} \mu_{\beta,\underline{\theta}}(B\cap [0,\frac{1}{\beta-1}) +i_{j})}{\sum_{j} \mu_{\beta,\underline{\theta}}([0,\frac{1}{\beta-1})+i_{j})}\),
so for \(\omega \in \{0,1\}^{n}\) ,
\[\bar{\mu}_{\beta,\underline{\theta}}[\omega] = \mathbf{1}M_{\omega}\begin{pmatrix}
        \mu_{\beta,\underline{\theta}}(F_{\beta,i}^{-1} ([0,\frac{1}{\beta-1}))+i_{0}) \\
    \vdots \\
    \mu_{\beta,\underline{\theta}}(F_{\beta,i}^{-1}([0,\frac{1}{\beta-1}))+i_{k})
\end{pmatrix} \frac{1}{\sum_{j} \mu_{\beta,\underline{\theta}}([0,\frac{1}{\beta-1})+i_{j})} .\]

For the sake of readability, let \(\mathbf{R} =\begin{pmatrix}
        \mu_{\beta,\underline{\theta}}(F_{i}^{-1}(B)+i_{0}) \\
    \vdots \\
    \mu_{\beta,\underline{\theta}}(F_{i}^{-1}(B)+i_{k})
\end{pmatrix} \frac{1}{\sum_{j} \mu_{\beta,\underline{\theta}}([0,\frac{1}{\beta-1})+i_{j})} \).

We define the \(n^{th}\) step potential of the measure \(\bar{\mu}_{\underline{\theta}}\) to be \[ \bar{\phi}_{n}(x) \defeq \log\left( \frac{\mathbf{1} M_{x_{1}}\dots M_{x_{n}} \mathbf{R}}{\mathbf{1} M_{x_{2}}\dots M_{x_{n}}\mathbf{R}}\right).\]
 
Note that there exists a positive vector \((1/2^{k+1},\dots,1/2^{k+1})\) such that \\ \( (1/2^{k+1},\dots,1/2^{k+1})\mathbf{1} = 1\). \\
This allows us to re express \(\bar{\phi}_{n}(x) = \log\left( \frac{\mathbf{1} M_{x_{1}}(1/2^{k+1},\dots,1/2^{k+1})\mathbf{1}M_{x_{1}}\dots M_{x_{n}} \mathbf{R}}{\mathbf{1} M_{x_{2}}\dots M_{x_{n}}\mathbf{R}}\right)\). 

Consider the variation of the \(n^{th}\) step potential of \(\bar{\mu}_{\underline{\theta}}\). 

\begin{align*}
    &\Bigg| \log \left( \frac{\mathbf{1} M_{x_{1}}(1/2^{k+1},\dots,1/2^{k+1})\mathbf{1}M_{x_{2}}\dots M_{x_{n}} \mathbf{R}}{\mathbf{1} M_{x_{2}}\dots M_{x_{n}}\mathbf{R}} \right) \\ & \qquad\qquad\qquad\qquad\qquad- \log \left(\frac{\mathbf{1} M_{x_{1}}(1/2^{k+1},\dots,1/2^{k+1})\mathbf{1}M_{x_{2}}\dots M_{x_{m}} \mathbf{R}}{\mathbf{1} M_{x_{2}}\dots M_{x_{m}}\mathbf{R}} \right) \Bigg| \\
    &\leq \Bigg| \log \left(  \frac{\mathbf{1} M_{x_{1}}(1/2^{k+1},\dots,1/2^{k+1})\mathbf{1}M_{x_{2}}\dots M_{x_{n}} \mathbf{R}}{\mathbf{1} M_{x_{2}}\dots M_{x_{n}}\mathbf{R}} \right) \Bigg| \\ &\qquad\qquad\qquad\qquad\qquad+ \Bigg| \log  \left(\frac{\mathbf{1} M_{x_{1}}(1/2^{k+1},\dots,1/2^{k+1})\mathbf{1}A_{x_{2}}\dots M_{x_{m}} \mathbf{R}}{\mathbf{1} A_{x_{2}}\dots M_{x_{m}}\mathbf{R}} \right)\Bigg| \\
    &=2\log\left( \mathbf{1} M_{x_{1}} (1/2^{k+1},\dots,1/2^{k+1}) \right)
\end{align*}

This implies \(\lvert \bar{\phi_{n}}(x) - \bar{\phi_{m}}(x)\rvert \leq 2\log(\lVert M_{x_{1}}\rVert ).\)

This allows us to establish the Weak Gibbs properties of \(\bar{\mu}_{\underline{\theta}}\) and so the local Gibbs properties \(\mu_{\underline{\theta}}\). 

\begin{lemma}\label{lemma : bargibbs}
    For \(\beta\in(1,2]\) a PV number, \(k\geq 1\), \(\theta = (f_{1},\dots,f_{k}), k > 1 \) and \(f_{i}\in \Q[\beta]\) for every \(i\). 
    Then \(\bar{\mu}_{\underline{\theta}}\) is the weak Gibbs measure associated to \(\bar{\phi}_{n}\).
    \begin{proof}
        Because \(\var(\bar{\phi_{n}}) \leq 2\log(\lVert M_{x_{1}}\rVert)\) and 
        \( \lvert \bar{\phi}(x) - \bar{\phi_{n}}(x) \rvert \leq \var{\phi_{n}}\) we see that the sequence \((n2\log(\lVert M_{x_{1}}\rVert))_{n}\) acts as a bound for the Weak Gibbs inequality as \(\lim_{n\to\infty}\frac{\log(n2\log(\lVert M_{x_{1}}\rVert))}{n} \to 0.\)
    \end{proof}
\end{lemma}

\begin{theorem}\label{theorem : locally weak gibbs}
    For \(\beta \in (1,2]\) a PV number, \(k\geq 1\) and \(\underline{\theta} = (f_{1},\dots,f_{k}), k>1\) with \(f_{i}\in\Q[\beta]\) for every \(i\) then \(\mu_{\underline{\theta}}\) is a locally Weak Gibbs measure associated to \(\phi\). 
    \begin{proof}
        By Lemma~\ref{lemma : bargibbs} we know that the measure \(\bar{\mu}_{\beta,\underline{\theta}}\) is Weak Gibbs. 
        By (Theorem~2.5~\cite{olivier}), whenever \(\bar{\mu}_{\beta,\underline{\theta}}\) is Weak Gibbs then \(\mu_{\beta,\underline{\theta}}\) is locally Weak Gibbs. 
    \end{proof}
\end{theorem}

\subsection{Pressure Result}
We now turn to establishing the key result of this work. 
We look at the pressure of the system \(P(l_{\underline{\theta}},T\vert_{l_{\underline{\theta}}},\phi)\) and relate this to the growth rate of the number of exact overlaps, \(\mathcal{N}\).
We begin by restating the definition of topological pressure Definition~\ref{definition : pressure definition}. 
For the space \(X_{\underline{\theta}}\) with the map \(\sigma : X_{\underline{\theta}} \to X_{\underline{\theta}}\) and the potential function \(\phi : X_{\underline{\theta}} \to \R \), define the pressure of \(\phi\) on \(X_{\underline{\theta}}\) under \(\sigma\) as 
\begin{equation*}
    P(X_{\underline{\theta}},\sigma,\phi) = \lim_{n \to \infty}\frac{1}{n} \log\left(\sum_{i_{1}...i_{n} \in X_{\underline{\theta}}}\exp \left(\sup_{\omega \in [i_{1}...i_{n}]}\sum_{j = 0}^{n-1}\phi(\sigma^{j}\omega)\right)\right).
    \end{equation*}

The use of restricting to \(X_{\underline{\theta}}\) is that \(B_{k+1}^{\N}\) has both recoverable and irrecoverable sequences combined in a single system. 
While, we know that \(l_{\underline{\theta}}\) contains the points which correspond with exact overlaps. 
Furthermore, we have the potential function \(\phi\) relates to the growth rate of the number of exact overlaps Lemma~\ref{lemma : sup n equiv}. 
These facts combined motivate the following key theorem. 

\begin{theorem}\label{theorem : pressure}
Let \(\beta\in(1,2]\) a PV number, \(k\geq 1\), \(\underline{\theta} = (f_{1},\dots,f_{k})\), \(f_{i}\in \Q[\beta]\) for all \(i\) then, 
the pressure function \(P(l_{\underline{\theta}},T\vert_{l_{\underline{\theta}}},\phi)\) satisfies
\[  \mathcal{N}\left(D,\{\frac{x+i}{\beta}  \colon i \in D\}\right) \leq (l_{\underline{\theta}},\sigma\vert_{l_{\underline{\theta}}},\phi).\]
    \begin{proof}
        Fix \(\beta \in (1,2]\) a PV number and let \(\underline{\theta} = (f_{1},\dots,f_{k})\) be a vector of \(\beta\) rational functions.
        
        Let \(\Omega_{n}\defeq [-\frac{f_{1}}{q_{1}}\beta^{n},\frac{f_{1}}{q_{1}}\beta^{n}]\times \cdots \times[-\frac{f_{k}}{q_{k}}\beta^{n},\frac{f_{k}}{q_{k}}\beta^{n}]\times[-\beta^{n},\beta^{n}]\ \subset \R^{k+1}.\)
        We identify \(\Omega_{n}\) with \(B_{k+1}^{n}\) for all \(n \in \N\) through rescaling by the vector \((\frac{f_{1}}{q_{k}},\dots,\frac{f_{k}}{q_{k}},1)\) and the projection map \(\pi^{-1}_{\beta}\).
        Considering the definition of recoverable pairs of words, we see that the region \( f_{1}x_{1}+\dots + f_{k}x_{k} - x_{k+1} \pm \frac{1}{\beta-1} = 0 \) bounds the region of \(\Omega_{n}\) in which recoverable pairs of words can exist.
        We re-scale \(\Omega_{n}\) to the cube \([-\frac{1}{\beta-1},\frac{1}{\beta-1}]^{k+1}\). 
        This scaling maps tubular neighbourhood of the plane, \( \frac{1}{q_{1}}x_{1}+\dots + \frac{1}{q_{k}}x_{k} - x_{k+1} \pm \frac{1}{\beta-1}= 0 \), to \( f_{1}x_{1}+\dots + f_{k}x_{k} - x_{k+1} \pm \frac{1}{\beta^{n}}\frac{1}{(\beta-1)^{2}} = 0 \). 
        Taking the limit as \(n\to \infty\) the tubular region, \( f_{1}x_{1}+\dots + f_{k}x_{k} - x_{k+1} \pm \frac{1}{\beta^{n}}\frac{1}{(\beta-1)^{2}} = 0 \), limits to \( \frac{1}{q_{1}}x_{1}+\dots + \frac{1}{q_{k}}x_{k} - x_{k+1} = 0 \).
        The plane \( \frac{1}{q_{1}}x_{1}+\dots + \frac{1}{q_{k}}x_{k} - x_{k+1} = 0 \) is a subset of \(l_{\underline{\theta}}.\) 
        Note that the closure of \( \frac{1}{q_{1}}x_{1}+\dots + \frac{1}{q_{k}}x_{k} - x_{k+1} = 0 \) with respect to the \(\beta\)ing map, \(T_{\beta}\) is \(l_{\underline{\theta}}.\)
        As all \(q_{i}\) are rational values this closure under the shift maps forms a hyperplane in the \(k+1\) torus, justifying the name \(\beta\) rational plane.
        By Corollary~\ref{corollary : lim sup equiv}, \[\log\frac{\mathcal{N}_{n}(D,D')}{\mathcal{N}_{n-1}(D,D')} \leq  \sum_{x\in D_{n}}\sup_{y\in [x]}\phi(y) ,\] and noting that \(\sum_{i=0}^{n} \log\frac{\mathcal{N}_{i}(D,D')}{\mathcal{N}_{i-1}(D,D')} = \log \mathcal{N}_{n}(D,D')\) we see
        \[\mathcal{N}_{n}(D,D') \leq \sum_{i_{1}...i_{n} \in X_{\underline{\theta}}}\exp \left(\sup_{\omega \in [i_{1}...i_{n}]}\sum_{j = 0}^{n-1}\phi(\sigma^{j}\omega)\right).\]
        Finally as,
        \(\mathcal{N}(D,\{x\to\frac{x+i}{\beta} \colon i \in D\}) = \lim_{n\to\infty}\frac{1}{n}\log\mathcal{N}(D,\{x\to\frac{x+i}{\beta} \colon i \in D) \implies \)
        \begin{align*}
        \mathcal{N}(D,\{x\to\frac{x+i}{\beta} \colon& i \in D\})\\ &\leq \lim_{n \to \infty}\frac{1}{n} \log\left(\sum_{i_{1}...i_{n} \in X_{\underline{\theta}}}\exp \left(\sup_{\omega \in [i_{1}...i_{n}]}\sum_{j = 0}^{n-1}\phi(\sigma^{j}\omega)\right)\right)\\ &=  P(l_{\underline{\theta}},T\vert_{l_{\underline{\theta}}},\phi).
        \end{align*}
    \end{proof}  
\end{theorem}

Combining the above theorem with Equation~\ref{equation : N above n} we get the following result for the dimension of the measure in terms of its pressure. 

\begin{theorem}\label{theorem: dim drop bound}
    Let \(\beta \in (1,2]\) be a PV number, \(k\geq 1\), \(\underline{\theta} =(f_{1},\dots,f_{k})\) with \(f_{i}\in \Q[\beta]\) for all \(i\).
    Then
    \[1 > \dim(\mu_{\beta,\underline{\theta}} )\geq \frac{2\log( k+2)}{\log\beta} - \frac{P(l_{\underline{\beta}},T\vert_{l_{\underline{\theta}}},\phi)}{\log \beta}>0. \]
    \begin{proof}
        By Theorem~\ref{lemma : sup n equiv}, we know \(\dim(\mu_{\beta,\theta})<1\). 
        Theorem~\ref{equation : N above n} concludes that  \(\frac{ H_{RW}(\underline{\theta})}{\log\beta}  \geq \frac{\mathcal{N}(D,\{\frac{x+i}{\beta} : i \in D\})}{\log\beta}\).
        Combining this with Theorem \ref{theorem : pressure}, the above and Hochman's result, (\(\dim(\mu_{\underline{\theta}}) = \min\{1,\frac{h_{RW}(\underline{\theta})}{\log\beta}\}\)).
        Gives,
        \[\log (k+2) \leq \frac{\mathcal{N}_{n}(D,D')}{\mathcal{N}_{n-1}(D,D')}\leq \sum_{x\in D_{n}} \sup_{y\in[x]}\phi(y),\]
        Which implies the following.
        \[P(l_{\underline{\theta}},T\vert_{l_{\underline{\theta}}},\phi)<\lim_{n\to\infty}\frac{1}{n}\log\left(\sum_{i_{1}\dots,i_{n}\in X_{\underline{\theta}}}(k+2)^{n}\right)\leq2\log(k+2)\]
        this gives the result.
    \end{proof}
\end{theorem}

We have now proven our key result that the dimension drop of \(\mu_{\underline{\theta}}\) is upper bounded by the pressure of \(\phi\) on the torus, restricted to varying planes \(l_{\underline{\theta}}\). 
This gives a characterisation of dimension for a wide range of self similar measures with overlaps in terms of the topological pressure of slices of a single dynamical system. 
It is the author's hope that this can be used to gain a deeper understanding of dimension drop and potentially a useful method for estimate dimension drop in these cases.
To this end we now include a comparison of the result with another method in the literature. 

\section{Dimension Drop for Padded Golden Mean}
The golden ratio or golden mean is a well known example of a PV number.
The golden mean is often denoted \(\varphi\) and we follow this convention to provide a clear distinction from an arbitrary PV number.

The golden mean  \(\varphi\) is the largest root of \(x^{2}-x-1=0\) and has the value of \(\frac{1+\sqrt{5}}{2}.\)
By rearranging the minimum polynomial it is relatively apparent that \(\pi_{\varphi}(100)=\pi_{\varphi}(011).\)
Due to these properties, amongst other key properties, it is one of the most well studied examples for numeration systems with a PV base. 

In this section we shall consider the IFS formed by the following three maps.
\[g_{0}(x) =\frac{x}{\varphi},\quad g_{\varphi}(x)=\frac{x+\frac{1}{\varphi}}{\varphi},\quad g_{1}(x)=\frac{x+1}{\varphi}.\]

As \(\Q[\beta]=\Q[\beta^{-1}]\) or by noting \(\frac{1}{\varphi} = \varphi-1\), the self similar measure \(\mu_{\varphi,\frac{1}{\varphi}}\) meets all the assumptions of this work. 

Now we consider two methods of bounding \(\dim(\mu_{\varphi,\frac{1}{\varphi}}).\)

As outlined in (Theorem 4.1 \cite{laulq}) there exists an algorithm to compute the \(L^{2}\) dimension of the measures \(\mu_{\beta,\underline{\theta}}.\)
This method involves building a transition matrix and taking the leading eigenvalue of its irreducible component. 
As noted after Theorem~\ref{lemma: invariant equiv} this transition matrix is a rescaling of the transition matrix for the automota constructed in Section~\ref{section: dynamic invariance}.

We now give the automaton for the system \(g_{0},g_{\varphi},g_{1}\) on the next page.

\begin{sidewaysfigure}
\centering

\resizebox{22cm}{!}{
\begin{tikzpicture}
\usetikzlibrary{positioning}

\tikzset{    
    every node/.append style = {
        draw = none,
        text = black,
    },
    every edge/.append style = {
        draw = black,
        arrows = {-Latex[scale = 1.5]}
    },
    every loop/.append style = {
        min distance= 10mm, 
        in=75, out=105, looseness=10,
        draw = black,
        arrows = {-Latex[scale = 1.5]}        
    },
    circ/.style = {
        draw = black,
        shape = circle,
        inner sep = 1pt,
        scale= .75,
        align=center,
        minimum size=1.2cm
    }
}

\node[circ] (10) at (0,0) {\(0\)};

\node[circ] (9) at (-3,12) {\(-2\varphi+3\)};
\node[circ] (8) at (-6,12) {\(+\varphi -2\)};
\node[circ] (7) at (-9,0) {\(-\varphi +1\)};
\node[circ] (6) at (-12,12) {\(2\varphi-4\)};
\node[circ] (5) at (-15,12) {\(-3\varphi +4\)};
\node[circ] (4) at (-18,0) {\(-1\)};
\node[circ] (3) at (-21,12) {\(-2\varphi +2\)} ;
\node[circ] (2) at (-24,12) {\(+\varphi -3\)} ; 
\node[circ] (1) at (-27,0) {\(-\varphi\)} ;

\node[circ] (11) at (3,12) {\(+\varphi -3\)};
\node[circ] (12) at (6,12) {\(-\varphi +2\)};
\node[circ] (13) at (9,0) {\(+\varphi -1\)};
\node[circ] (14) at (12,12) {\(-2\varphi+4\)};
\node[circ] (15) at (15,12) {\(+3\varphi-4\)};
\node[circ] (16) at (18,0) {\(+1\)};
\node[circ] (17) at (21,12) {\(+2\varphi-2\)};
\node[circ] (18) at (24,12) {\(-\varphi+3\)};
\node[circ] (19) at (27,0) {\(+\varphi\)};

\draw (1) edge [loop,edge label =+1] (1) ;

\draw (2) edge [edge label = \(+\varphi-1\)] (1) ;
\draw (2) edge [bend right, edge label = \(+1\)] (3) ;

\draw (3) edge [bend right, edge label = \(+\varphi-1\)] (2) ;
\draw (3) edge [bend right, edge label = \(+1\)] (4) ;

\draw (4) edge [edge label = \(0\)] (1) ;
\draw (4) edge [edge label = \(-\varphi +2\) ] (3) ;
\draw (4) edge [loop, edge label = \(+\varphi-1\) ] (4) ;
\draw (4) edge [edge label = \(+1\) ] (7) ;

\draw (5) edge [bend right, edge label = \(0\)] (2) ;
\draw (5) edge [edge label = \(-\varphi +2\)] (4) ;
\draw (5) edge [bend right, edge label = \(+\varphi-1\)] (6) ;
\draw (5) edge [bend right, edge label = \(+1\)] (8) ;

\draw (6) edge [edge label = \(+\varphi-2\)] (1) ;
\draw (6) edge [bend right, edge label = \(0\)] (3) ;
\draw (6) edge [bend right, edge label = \(-\varphi+2\)] (5) ;
\draw (6) edge [edge label = \(+\varphi-1\)] (7) ;
\draw (6) edge [bend right, edge label = \(+1\)] (9) ;

\draw (7) edge [bend left, edge label = \(-\varphi+1\)] (1) ;
\draw (7) edge [edge label = \(+\varphi -2\)] (2) ;
\draw (7) edge [bend right, edge label = \(0\)] (4) ;
\draw (7) edge [loop, edge label = \(-\varphi+2\)] (7) ;
\draw (7) edge [bend right, edge label = \(+\varphi-1\)] (8) ;
\draw (7) edge [bend left, edge label = \(+1\)] (10) ;

\draw (8) edge [edge label = \(-1\) ] (1) ;
\draw (8) edge [bend right, edge label = \(-\varphi+1\) ] (3) ;
\draw (8) edge [edge label = \(+\varphi-2\) ] (4) ;
\draw (8) edge [bend right, edge label = \(0\) ] (7) ;
\draw (8) edge [bend right, edge label = \(-\varphi+2\) ] (9) ;
\draw (8) edge [edge label = \(+\varphi-1\) ] (10) ;
\draw (8) edge [bend right, edge label = \(+1\) ] (12) ;

\draw (9) edge [bend right, edge label = \(-1\) ] (2) ;
\draw (9) edge [bend right, edge label = \(-\varphi+1\) ] (4) ;
\draw (9) edge [bend right, edge label = \(+\varphi-2\) ] (6) ;
\draw (9) edge [bend right, edge label = \(0\) ] (8) ;
\draw (9) edge [edge label = \(-\varphi+2\) ] (10) ;
\draw (9) edge [bend right, edge label = \(+\varphi-1\) ] (11) ;
\draw (9) edge [, edge label = \(+1\) ] (13) ;

\draw (10) edge [bend left, edge label = \(-1\) ] (4) ;
\draw (10) edge [edge label = \(-\varphi+1\) ] (7) ;
\draw (10) edge [bend left, edge label = \(+\varphi-2\) ] (8) ;
\draw (10) edge [loop, edge label = \(0\) ] (10) ;
\draw (10) edge [bend right, edge label = \(-\varphi+2\) ] (12) ;
\draw (10) edge [edge label = \(+\varphi-1\) ] (13) ;
\draw (10) edge [bend right, edge label = \(+1\) ] (16) ;

\draw (11) edge [bend right, edge label = \(-1\) ] (8) ; 
\draw (11) edge [bend right, edge label = \(-\varphi+1\) ] (9) ;
\draw (11) edge [edge label = \(+\varphi-2\) ] (10) ;
\draw (11) edge [bend right, edge label = \(0\) ] (12) ;
\draw (11) edge [bend right, edge label = \(-\varphi+2\) ] (14) ;
\draw (11) edge [edge label = \(+\varphi-1\) ] (16) ;
\draw (11) edge [bend left, edge label = \(+1\) ] (18) ;

\draw (12) edge [bend right, edge label = \(-1\) ] (8) ;
\draw (12) edge [edge label = \(-\varphi+1\) ] (10) ;
\draw (12) edge [bend right, edge label = \(+\varphi-2\) ] (11) ;
\draw (12) edge [edge label = \(0\) ] (13) ;
\draw (12) edge [, edge label = \(-\varphi+2\) ] (16) ;
\draw (12) edge [bend left, edge label = \(+\varphi-1\) ] (17) ;
\draw (12) edge [edge label = \(+1\) ] (19) ;

\draw (13) edge [bend right, edge label = \(-1\) ] (10) ;
\draw (13) edge [bend left, edge label = \(-\varphi+1\) ] (12) ;
\draw (13) edge [loop, edge label = \(+\varphi-2\) ] (13) ;
\draw (13) edge [bend left, edge label = \(0\) ] (16) ;
\draw (13) edge [edge label = \(-\varphi+2\) ] (18) ;
\draw (13) edge [bend right, edge label = \(+\varphi-1\) ] (19) ;

\draw (14) edge [bend right, edge label = \(-1\) ] (11) ;
\draw (14) edge [edge label = \(-\varphi+1\) ] (13) ;
\draw (14) edge [bend left, edge label = \(+\varphi-2\) ] (15) ;
\draw (14) edge [bend left, edge label = \(0\) ] (17) ;
\draw (14) edge [edge label = \(-\varphi+2\) ] (19) ;

\draw (15) edge [bend right, edge label = \(-1\) ] (12) ;
\draw (15) edge [bend left, edge label = \(-\varphi+1\) ] (14) ;
\draw (15) edge [edge label = \(+\varphi-2\) ] (16) ;
\draw (15) edge [bend left, edge label = \(0\) ] (18) ;

\draw (16) edge [edge label = \(-1\) ] (13) ;
\draw (16) edge [loop, edge label = \(-\varphi+1\) ] (16) ;
\draw (16) edge [edge label = \(+\varphi-2\) ] (17) ;
\draw (16) edge [edge label = \(0\) ] (19) ;

\draw (17) edge [bend left, edge label = \(-1\) ] (16) ;
\draw (17) edge [bend right, edge label = \(-\varphi+1\) ] (18) ;

\draw (18) edge [edge label = \(-1\) ] (19) ;
\draw (18) edge [bend right, edge label = \(-\varphi+1\) ] (17) ;

\draw (19) edge [loop , edge label = \(-1\) ] (19);

\end{tikzpicture}
}
\end{sidewaysfigure}

The transition matrix associated to that automaton is,  
\[
    \begin{bmatrix}
    0 & 1 & 0 & 0 & 0 & 0 & 0 & 0 & 0 & 0 & 0 & 0 & 0 & 0 & 0 & 0 & 0 \\
    1 & 0 & 1 & 0 & 0 & 0 & 0 & 0 & 0 & 0 & 0 & 0 & 0 & 0 & 0 & 0 & 0 \\
    0 & 1 & 1 & 0 & 0 & 1 & 0 & 0 & 0 & 0 & 0 & 0 & 0 & 0 & 0 & 0 & 0 \\
    3 & 0 & 1 & 0 & 1 & 0 & 1 & 0 & 0 & 0 & 0 & 0 & 0 & 0 & 0 & 0 & 0 \\
    0 & 3 & 0 & 1 & 0 & 1 & 0 & 1 & 0 & 0 & 0 & 0 & 0 & 0 & 0 & 0 & 0 \\
    1 & 0 & 3 & 0 & 0 & 1 & 1 & 0 & 1 & 0 & 0 & 0 & 0 & 0 & 0 & 0 & 0 \\
    0 & 1 & 0 & 0 & 0 & 3 & 0 & 1 & 1 & 0 & 1 & 0 & 0 & 0 & 0 & 0 & 0 \\
    1 & 0 & 1 & 0 & 1 & 0 & 3 & 0 & 1 & 1 & 0 & 1 & 0 & 0 & 0 & 0 & 0 \\
    0 & 0 & 1 & 0 & 0 & 1 & 1 & 0 & 3 & 1 & 1 & 1 & 0 & 0 & 1 & 0 & 0 \\
    0 & 0 & 0 & 0 & 0 & 1 & 0 & 1 & 1 & 0 & 3 & 0 & 1 & 0 & 1 & 0 & 1 \\
    0 & 0 & 0 & 0 & 0 & 0 & 1 & 0 & 1 & 0 & 0 & 3 & 0 & 0 & 1 & 1 & 0 \\
    0 & 0 & 0 & 0 & 0 & 0 & 0 & 0 & 1 & 0 & 1 & 1 & 0 & 0 & 3 & 0 & 1 \\
    0 & 0 & 0 & 0 & 0 & 0 & 0 & 0 & 0 & 1 & 0 & 1 & 0 & 1 & 0 & 3 & 0 \\
    0 & 0 & 0 & 0 & 0 & 0 & 0 & 0 & 0 & 0 & 1 & 0 & 1 & 0 & 1 & 0 & 3 \\
    0 & 0 & 0 & 0 & 0 & 0 & 0 & 0 & 0 & 0 & 0 & 1 & 0 & 0 & 1 & 1 & 0 \\
    0 & 0 & 0 & 0 & 0 & 0 & 0 & 0 & 0 & 0 & 0 & 0 & 0 & 0 & 1 & 0 & 1 \\
    0 & 0 & 0 & 0 & 0 & 0 & 0 & 0 & 0 & 0 & 0 & 0 & 0 & 0 & 0 & 1 & 0
    \end{bmatrix}
\]

We now rescale such that the transition matrix is probabilistic in its maximum row.
Then the leading eigenvalue of the rescaled matrix can be calculated and is \(\lambda =  0.6343.\)
This leads to that \(D_{\mu_{\beta,\frac{1}{\varphi}}}= \frac{\log (0.6343)}{\log(\frac{1}{\varphi})}=0.946014\).

Now we use Theorem \ref{theorem: dim drop bound} to calculate a lower bound of \(\dim(\mu_{\varphi,\frac{1}{\varphi}})\) instead.

As noted prior to this section a complete general analysis of this potential function restricted to any rational slice could yield insight into dimension drop for this entire class of measures. 
In lieu of this general framework we must deal with the measures on a cases by case basis. 

For \(\mu_{\varphi,\frac{1}{\varphi}}\) we note that for the \(\varphi\) rational function \(f_{\varphi}=\frac{1}{\varphi}\) has \(q_{\varphi}=1\). 
As such \(l_{\frac{1}{\varphi}}\) is the line \(x=y \subset [-\varphi,\varphi]^{2}.\)

Considering the symbolic encoding of \([-\varphi,\varphi]^{2}\) into \(\{-1,0,1\}^{\N}\times\{-10,1\}^{\N}\) given by \(\pi_{\varphi}.\) 
Let \(a =\begin{pmatrix}
    0 \\0
\end{pmatrix} \),\(b = \begin{pmatrix}
    1 \\1
\end{pmatrix}\) and \(c=\begin{pmatrix}
    -1 \\-1
\end{pmatrix}.\)
Then we note that \(l_{\frac{1}{\varphi}} \subset[\omega]\) for \(\omega \in \{a,b,c\}^{n}\) for all \(n\in \N\).

Therefore it suffices to consider terms of the form 
\[\phi\vert_{n}(\underline{x}) = \log\frac{(1,1,1,1)A_{x_{1}}\dots A_{x_{n}}(1,0,0,0)^{T}}{(1,1,1,1)A_{x_{2}}\dots A_{x_{n}}(1,0,0,0)^{T}},\]
for \(\underline{x}\in \{a,b,c\}^{n}.\) 

This allows for the use a computer to calculate \(\phi\) for this collection of cylinders.  

From this we are able to conclude that \(P(X_{\frac{1}{\varphi},\sigma,\phi}) \leq \frac{\log7}{\log\phi}.\)
Therefore \[\dim(\mu_{\varphi,\frac{1}{\varphi}})\geq \frac{2\log 3}{\log\varphi}-\frac{\log7}{\log\varphi}=0.5222 \]

While \(0.5222 < 0.946014,\) we note that the calculation of matrix and its eigenvalue are bespoke for each PV number as such existing methods do not allow the following corollary to be drawn.

\begin{corollary}
    For \(\beta\in(1,2)\) a PV number then,
    \[\dim(\mu_{\beta,\frac{1}{\beta}})\geq  \frac{2\log 3}{\log\beta}-\frac{\log7.5}{\log\beta}.\]
\end{corollary}
This is as for any \(f_{i}\) with \(q_{i}=1\) the above calculation holds completely, in fact a better could be calculated using these methods but again would rely on a deeper understanding of the products of the matrices \(A_{i}.\)

\bibliographystyle{plain} 
\bibliography{ref_frac}

\end{document}